# ERROR BOUNDS FOR CONTROL CONSTRAINED SINGULARLY PERTURBED LINEAR-QUADRATIC OPTIMAL CONTROL PROBLEMS*

SEI HOWE† AND PANOS PARPAS‡

**Abstract.** We present a methodology for bounding the error term of an asymptotic solution to a singularly perturbed optimal control (SPOC) problem whose exact solution is known to be computationally intractable. In previous works, reduced or computationally tractable problems that are no longer dependent on the singular perturbation parameter $\epsilon$, where $\epsilon$ represents a small, non-negative number, have provided asymptotic error bounds of the form $O(\epsilon)$. Specifically, the optimal solution $\bar{V}$ of the reduced problem has been shown to be asymptotically equivalent in $\epsilon$ to the optimal solution $V(\epsilon)$ of the singularly perturbed problem in the sense that $|V(\epsilon) - \bar{V}| = O(\epsilon)$ as $\epsilon \to 0$. In this paper, we improve on this result by incorporating a duality theory into the SPOC problem and derive an upper bound $\chi_u(\epsilon)$ and lower bound $\chi_l(\epsilon)$ of $V(\epsilon)$ that hold for arbitrary $\epsilon$ and, furthermore, satisfy the inequality $|\chi_u(\epsilon) - \chi_l(\epsilon)| \leq C\epsilon$ for small $\epsilon$, with the constant $C$ determined. We carry out numerical experiments to illustrate the computational savings obtained for the upper and lower bound. In particular, we generate a set of 50 random SPOC problems of a specific form and show that for $\epsilon$ smaller than $10^{-2}$, it becomes faster, on average, to solve for the bounds rather than the SPOC problem and for $\epsilon = 10^{-5}$, the computational time for the upper and lower bounds is approximately 20 times faster, on average, than that of the SPOC problem.

## 1. Introduction.

Singularly perturbed optimal control (SPOC) problems are known to be computationally challenging to solve; accordingly, various techniques have been developed in order to obtain approximations to the optimal solution. In particular, the idea of solving a computationally tractable and dimensionally reduced problem has been heavily investigated. For a wide class of control constrained SPOC problems with both linear and non linear dynamics, a reduced problem has been derived whose optimal solution $\bar{V}$ is shown to be the limit of the optimal solution $V(\epsilon)$ of the SPOC problem (see [2] - [5], [10], [15], [27]); specifically

$$\lim_{\epsilon \to 0} V(\epsilon) = \bar{V}. \tag{1}$$

In this paper, we improve on the result in (1) for a special class of linear control constrained SPOC problems by determining an upper bound $\chi_u(\epsilon)$ and a lower bound $\chi_l(\epsilon)$ on $V(\epsilon)$ that hold for arbitrary values of $\epsilon$ and, futhermore, satisfy the condition

$$|\chi_l(\epsilon) - \chi_u(\epsilon)| \leq C\epsilon \quad \text{for all } \epsilon \in (0, \epsilon^*], \tag{2}$$

where the constant $C$ is explicitly determined.

As the SPOC problem we consider is a linear-quadratic minimisation problem, an upper bound is easily found by evaluating the original problem with the optimal control of an appropriately constructed reduced problem provided that it is a feasible control. However, there has been no comparable result for the lower bound and no ensuing result for the constant in (2). Deriving this lower bound and constant has been difficult, in part, due to the lack of a duality framework in which to formulate

*The authors acknowledge support from EPSRC grants EP/M028240, EP/K040723 and an FP7 Marie Curie Career Integration Grant (PCIG11-GA-2012-321698 SOC-MP-ES)

†Department of Computer Science, Imperial College London, U.K. (sei.howe11@imperial.ac.uk). Questions, comments, or corrections to this document may be directed to that email address.

‡Department of Computer Science, Imperial College London, U.K.(p.parpas@imperial.ac.uk)

the SPOC problem. By extending the duality techniques in [1] and [7] to the case of SPOC problems, we develop such a framework and derive an explicit dual problem with the strong duality property. As the dual problem is a maximisation problem, it follows from strong duality that any feasible control will provide a lower bound on the optimal solution of the SPOC problem. In order to show that the error bound in (2) holds, we will use the strong duality result and the reduced problem to contsruct a control that provides an $O(\epsilon)$ approximation to the optimal control of the dual problem as $\epsilon \to 0$ uniformly on the entire time interval. By evaluating the SPOC primal and dual problems with the optimal control of the reduced problem and the constructed control respectively, one obtains the error bound on $V(\epsilon)$ of the form (2) with the constant $C$ determined.

The novel aspects of this paper are the bound on the error term of an asymptotic solution to a control constrained SPOC problem using computationally feasible reduced solutions and the introduction of a dual optimal control problem satisfying the strong duality property. Using our framework, it becomes possible to obtain error bounds on the SPOC problem for any arbitrary value of $\epsilon$ as well as a criteria for evaluating how good of an approximation the reduced problem provides. We also show that significant computational savings may be achieved when obtaining the upper and lower bounds instead of solving the primal or dual SPOC problem. Specifically, for a set of 50 randomly generated SPOC problems, when $\epsilon \leq 10^{-2}$, the bounds start to be faster to compute than either the primal or dual SPOC problem. Furthermore, for $\epsilon = 10^{-5}$, the the computational time for the bounds is approximately 160 and 20 times faster than that of solution to the primal and dual SPOC problems respectively.

**1.1. Formulation of the primal problem.** We consider the following SPOC problem

$$(\mathbf{P}) \qquad \begin{cases} \underset{\hat{u}\in U}{\text{minimise}} & J^{\mathbf{P}}(\hat{z}_1, \hat{z}_2, \hat{u}, \epsilon), \\ \text{subject to} & \dfrac{d\hat{z}_1}{dt} = A_{11}(t,\epsilon)\hat{z}_1(t,\epsilon) + A_{12}(t,\epsilon)\hat{z}_2(t,\epsilon) + b_1(t,\epsilon)\hat{u}(t,\epsilon), \\ & \epsilon\dfrac{d\hat{z}_2}{dt} = A_{21}(t,\epsilon)\hat{z}_1(t,\epsilon) + A_{22}(t,\epsilon)\hat{z}_2(t,\epsilon) + \epsilon b_2(t,\epsilon)\hat{u}(t,\epsilon), \\ & \hat{z}_1(0,\epsilon) = z_{1,0}, \quad \hat{z}_2(0,\epsilon) = z_{2,0}, \quad \text{for } \epsilon \in (0,\epsilon^*], \end{cases}$$

where the functional, $J^{\mathbf{P}}$, is defined as

$$J^{\mathbf{P}} = \frac{1}{2}\int_0^1 \hat{z}(t,\epsilon)^T Q(t,\epsilon)\hat{z}(t,\epsilon) + \hat{u}(t,\epsilon)^T R(t,\epsilon)\hat{u}(t,\epsilon)dt + \frac{1}{2}\hat{z}(1,\epsilon)^T\pi(\epsilon)\hat{z}(1,\epsilon). \tag{3}$$

and $\hat{z}^T = \left[\hat{z}_1^T, \hat{z}_2^T\right]$. We let $V^{\mathbf{P}}(\epsilon)$ denote the value of $J^{\mathbf{P}}$ evaluated at the optimal control, denoted by $u$, for fixed $\epsilon$. For all $\epsilon \in (0,\epsilon^*]$, $\hat{z}_1 \in W^{1,2}([0,1];\mathbb{R}^m)$, $\hat{z}_2 \in W^{1,2}([0,1];\mathbb{R}^n)$ as a function of $t$, where the space $W^{1,2}$ denotes the Sobelov space of absolutely continuous functions. Note that for $\epsilon = 0$, the dimension of the problem $\mathbf{P}$ drops from $m+n$ to $m$ and the boundary condition for $\hat{z}_2$ may no longer be satisfied. For any $\epsilon \in (0,\epsilon^*]$, the control $\hat{u}$ is constrained to be in the set $U$, where

$$U = \{u \in L^2([0,1];\mathbb{R}^k) : \alpha_j(t) \leq u_j(t) \leq \beta_j(t) \text{ for all } j = 1,\dots k, t \in [0,1]\} \tag{4}$$

where the functions $\alpha, \beta$ are smooth for $t \in [0,1]$. We make some standard assumptions for linear-quadratic optimal control models, namely that $R$ is positive definite and $Q$ and $\pi$ are positive semi-definite, and further impose the conditions that $R$ is diagonal and that $Q$ and $\pi$ are symmetric and positive definite. The diagonality assumption is needed in order to formulate the dual objective function, and the symmetry and positive definiteness assumptions are used in the proof of strong duality. Futhermore, for $\epsilon \in [0, \epsilon^*]$, we assume that $\pi$ has the following block-diagonal structure

$$\pi(\epsilon) = \begin{bmatrix} \pi_{11}(\epsilon) & 0 \\ 0 & \epsilon\pi_{22}(\epsilon) \end{bmatrix}, \qquad \pi_{11} \in \mathbb{R}^{m\times m},\ \pi_{22} \in \mathbb{R}^{n\times n}.$$

This assumption on the structure of $\pi$ is helpful in order to easily obtain the reduced problem. Finally, we impose the following standard assumptions for singularly perturbed problems (see [10], [17] - [24]):

(a) the matrix $A_{22}(t, \epsilon)$ is negative definite for all $t \in [0, 1]$, $\epsilon \in [0, \epsilon^*]$,

(b) for any fixed $\epsilon \in [0, \epsilon^*]$, the matrices $Q$, $R$, $A_{i,j}$, and $b_i$ for $i, j = 1, 2$, are smooth for $t \in [0, 1]$,

(c) $Q$, $R$, $\pi$, $A_{ij}$, and $b_i$, for $i, j = 1, 2$, all have a convergent power series expansion in $\epsilon$ which is valid over their respective domains.

The feasible set for $\mathbf{P}$ can be written as

$$\Sigma = \Big\{ (\hat{z}, \hat{u}) : \hat{z} \in W^{1,2}([0, 1]; \mathbb{R}^{m+n}), \hat{u} \in U, \hat{z}(0, \epsilon) = z_0 \text{ for all } \epsilon, \tag{5}$$
$$I^\epsilon \frac{d\hat{z}}{dt} = A(t, \epsilon)\hat{z}(t, \epsilon) + I^\epsilon b(t, \epsilon)\hat{u}(t, \epsilon),\ \ t \in [0, 1],\ \ \epsilon \in (0, \epsilon^*] \Big\},$$

where $z_0^T = [z_{1,0}^T, z_{2,0}^T]$ and

$$A(t, \epsilon) = \begin{bmatrix} A_{11}(t, \epsilon) & A_{12}(t, \epsilon) \\ A_{21}(t, \epsilon) & A_{22}(t, \epsilon) \end{bmatrix}, \qquad A_{11} \in \mathbb{R}^{m\times m},\ A_{22} \in \mathbb{R}^{n\times n},$$
$$b(t, \epsilon) = \begin{bmatrix} b_1(t, \epsilon) \\ b_2(t, \epsilon) \end{bmatrix}, \qquad b_1 \in \mathbb{R}^{m\times k},\ b_2 \in \mathbb{R}^{n\times k},$$
$$I^\epsilon = \begin{bmatrix} I_m & 0 \\ 0 & \epsilon I_n \end{bmatrix}.$$

and $I_j$ is the $j \times j$ identity matrix for $j = m, n$.

Remark 1. *We assume the set $\Sigma$ is nonempty. Since the solution set $\Sigma$ is closed and convex and $J^{\mathbf{P}}$ is strictly convex, continuous and coercive over $\Sigma$, there exists a unique solution to the minimisation problem $\mathbf{P}$ (see [11], Chap. 2, Proposition 1.2]).*

Remark 2. *One may also consider the case where the initial conditions in $\mathbf{P}$ are functions of $\epsilon$ with only a few minor changes to the analysis in this paper. For simplicity, however, we have chosen initial conditions that are independent of $\epsilon$.*

Several authors have obtained reduced models with asymptotically equivalent optimal solutions to variations of the above problem. Some of the first reduced models with such a property were developed for the linear-quadratic, unconstrained control case and were obtained using the boundary layer method (see [17] - [18], [21] - [24]). More recently, various other methods have been developed in order to obtain reduced models for general, nonlinear, constrained control, SPOC problems (see [2] - [5], [10], [15], [27]). While a reduction method has been developed for a more complex problem than ours, the reasons for looking at $\mathbf{P}$ are three-fold:

1. We aim to formulate an explicit dual SPOC problem with the strong duality property and closed form expressions have only be obtained for the linear-quadratic control constrained case [1], [7].

2. Linear-quadratic, control constrained, SPOC problems are widely encountered in many areas such as aerospace engineering [16], [12], electrical circuits [4], [18], chemical processes [18] and many others.

3. Our method of obtaining an error bound of the form (2) requires closed form expressions for the optimality conditions. A problem of the form $\mathbf{P}$ can provide such optimality conditions.

Our arguments rely on adaptations of the boundary layer method to the case of constrained control and adaptations of a duality framework to SPOC problems; each of which have been made explicit for the first time, to the extent of the authors knowledge, in this paper.

**2. Statement of main results.** The statement of the results has been divided into two sections: the first covers the formulations of the associated optimal control problems, namely the dual and reduced problems, the second covers the theorems relating to the derivation of constant in the error bound in (2).

**2.1. Formulation of the dual and reduced problem.** We begin by presenting the formulation of the dual maximisation problem to $\mathbf{P}$ and then follow with the formulation of the reduced problem. The construction of the dual problem, based on that of the unperturbed case [1] and [7], is found in Section 5. The dual problem can be formulated as

$$(\mathbf{D}) \quad \begin{cases} \underset{\hat\rho_1,\hat\rho_2,\hat\gamma_1,\hat\gamma_2}{\text{maximise}} & J^{\mathbf{D}}(\hat\rho_1,\hat\rho_2,\hat\gamma_1,\hat\gamma_2,\epsilon), \\ \text{subject to} & \dfrac{d\hat\gamma_1}{dt} = -A_{11}(t,\epsilon)^T\hat\gamma_1(t,\epsilon) - A_{21}(t,\epsilon)^T\hat\gamma_2(t,\epsilon) + \hat\rho_1(t,\epsilon), \\ & \epsilon\dfrac{d\hat\gamma_2}{dt} = -A_{12}(t,\epsilon)^T\hat\gamma_1(t,\epsilon) - A_{22}(t,\epsilon)^T\hat\gamma_2(t,\epsilon) + \hat\rho_2(t,\epsilon), \end{cases}$$

where the functional $J^{\mathbf{D}}$ is given by

$$J^{\mathbf{D}} = \int_0^1 -\frac{1}{2}\hat\rho^T Q^{-1}\hat\rho - \theta(\hat\gamma)dt - \hat\gamma(0,\epsilon)^T I^\epsilon z_0 - \frac{1}{2}\hat\gamma(1,\epsilon)^T I^\epsilon \pi(\epsilon)^{-1} I^\epsilon \hat\gamma(1,\epsilon), \tag{6}$$

and $\hat\rho^T = \left[\hat\rho_1^T, \hat\rho_2^T\right]$, $\hat\gamma^T = \left[\hat\gamma_1^T, \hat\gamma_2^T\right]$. We let $V^{\mathbf{D}}(\epsilon)$ denote the value of $J^{\mathbf{D}}$ evaluated at the optimal control, denoted by $\rho$, for fixed $\epsilon$. In (6), the function $\theta(\hat\gamma) = \sum_{j=1}^k \theta_j(\hat\gamma)$ and $\theta_j$ is defined as

$$\theta_j(\hat\gamma) = \begin{cases} \frac{1}{2R_{jj}}(\hat\gamma^T I^\epsilon b)_j^2 & \text{if } \alpha_j \le (R^{-1}b^T I^\epsilon\hat\gamma)_j \le \beta_j, \\ \alpha_j(\hat\gamma^T I^\epsilon b)_j - \frac{1}{2}\alpha_j^2 R_{jj} & \text{if } \alpha_j > (R^{-1}b^T I^\epsilon\hat\gamma)_j, \\ \beta_j(\hat\gamma^T I^\epsilon b)_j - \frac{1}{2}\beta_j^2 R_{jj} & \text{if } \beta_j < (R^{-1}b^T I^\epsilon\hat\gamma)_j. \end{cases} \tag{7}$$

Note that we have omitted the dependence on $t$ and $\epsilon$ in the integrand of (6) and equation (7) for notational simplicity. For all $\epsilon \in (0,\epsilon^*]$, $\hat\gamma_1, \hat\rho_1 \in W^{1,2}([0,1];\mathbb{R}^m)$, and $\hat{\gamma}_2, \hat{\rho}_2 \in W^{1,2}([0,1];\mathbb{R}^n)$ as functions of $t$. The feasible set for $\mathbf{D}$ can be written as

$$\Sigma_1 = \left\{(\hat{\gamma},\hat{\rho}) : \hat{\gamma},\hat{\rho} \in W^{1,2}([0,1];\mathbb{R}^{m+n}), I^{\epsilon}\frac{d\hat{\gamma}}{dt} = -A(t,\epsilon)^T\hat{\gamma} + \hat{\rho}, t \in [0,1], \epsilon \in (0,\epsilon^*]\right\}.$$

Remark 3. *The objective functional $J^{\mathbf{D}}$ is not necessarily coercive over the feasible set $\Sigma_1$; hence, we cannot immediately conclude that a unique solution exists. However, as there exists a unique solution to $\mathbf{P}$ by remark 1, the strong duality result in Section 2.2 will lead to the existence of a unique solution for $\mathbf{D}$.*

The dual problem $\mathbf{D}$ and the strong duality result are instrumental in obtaining the constant in the error bound (2) as they imply $V^{\mathbf{P}}(\epsilon) = V^{\mathbf{D}}(\epsilon)$. Hence any lower bound for $V^{\mathbf{D}}(\epsilon)$ will also be a lower bound for $V^{\mathbf{P}}(\epsilon)$. We now formulate the reduced problem for $\mathbf{P}$. Consider the following problem

$$(\bar{\mathbf{P}}) \qquad \begin{cases} \underset{\hat{u}\in U}{\text{minimise}} & \bar{J}^{\mathbf{P}}(\hat{x},\hat{u}), \\ \text{subject to} & \frac{d\hat{x}}{dt} = \mathcal{A}(t)\hat{x}(t) + b_1^0(t)\hat{u}(t), \\ & \hat{x}(0) = z_{1,0}, \end{cases}$$

where the functional $\bar{J}^{\mathbf{P}}(\hat{x},\hat{u})$ is defined by

$$(8) \qquad \bar{J}^{\mathbf{P}}(\hat{x},\hat{u}) = \frac{1}{2}\int_0^1 \hat{x}(t)^T\mathcal{Q}(t)\hat{x}(t) + \hat{u}(t)^TR^0(t)\hat{u}(t)dt + \frac{1}{2}\hat{x}(1)^T\pi_{11}^0\hat{x}(1).$$

We let $\bar{V}^{\mathbf{P}}$ denote the value of $\bar{J}^{\mathbf{P}}$ evaluated at the optimal control, denoted by $\bar{u}$. The notation $\pi_{11}^0$ is used to denote the first term in the power series expansion for $\pi_{11}$. Similar notation holds for the first term in the expansion for all other matrices satisfying assumption (c). Note that $\hat{x} \in W^{1,2}([0,1];\mathbb{R}^m)$ as a function of $t$ and $\hat{u}$ is constrained to be in the set $U$ defined in (4). The matrices $\mathcal{Q}$ and $\mathcal{A}$ are defined as

$$(9) \qquad \begin{aligned} \mathcal{Q} &= Q_{11}^0 - Q_{12}^0(A_{22}^0)^{-1}A_{21}^0 - (A_{21}^0)^T((A_{22}^0)^{-1})^TQ_{21}^0 + (A_{21}^0)^T((A_{22}^0)^{-1})^TQ_{22}^0(A_{22}^0)^{-1}A_{21}^0, \\ \mathcal{A} &= A_{11}^0 - A_{12}^0(A_{22}^0)^{-1}A_{21}^0. \end{aligned}$$

For simplicity we have forgone the dependence of the variables in (9) on $t$. We assume the matrix $\mathcal{Q}$ in (9) is positive semi-definite. The feasible set for $\bar{\mathbf{P}}$ is defined as follows

$$(10) \qquad \bar{\Sigma} = \left\{(\hat{x},\hat{u}) : \hat{x} \in W^{1,2}([0,1];\mathbb{R}^m), \;\; \hat{u} \in U, \; \hat{x}(0) = z_{1,0} \right. \\ \left. \frac{d\hat{x}}{dt} = \mathcal{A}(t)\hat{x}(t) + b_1^0(t)\hat{u}(t), \;\; t \in [0,1]\right\}.$$

Remark 4. *As the feasible set of $\bar{\mathbf{P}}$ given in (10) is closed and convex and $\bar{J}^{\mathbf{P}}$ in (8) is strictly convex, continuous and coercive over the feasible set, there exists a unique solution to $\bar{\mathbf{P}}$ (see [11], Chap. 2, Proposition 1.2]).*

The problem $\bar{\mathbf{P}}$ is equivalent to $\mathbf{P}$ evaluated at $\epsilon = 0$ and with the boundary condition for $\hat{z}_2$ ignored. We regard the problem $\bar{\mathbf{P}}$ as the reduced problem associated with $\mathbf{P}$ and, in Section 3, we show that the necessary optimality conditions for $\mathbf{P}$ converge to those of $\bar{\mathbf{P}}$ as $\epsilon \to 0$.

**2.2. Main results.** The methodology for obtaining the asymptotic error bound of $\mathbf{P}$ is depicted in Fig. 1. Strong duality will provide the equalities $V^{\mathbf{P}}(\epsilon) = V^{\mathbf{D}}(\epsilon)$. This, along with the evaluation of the functionals $J^{\mathbf{P}}$ and $J^{\mathbf{D}}$ with the optimal control of $\bar{V}^{\mathbf{P}}$, denoted by $\bar{u}$, and a constructed control $\boldsymbol{\rho}(\bar{u})$ respectively, provides the error bound of the form (2) with the constant determined. The convergence of the optimality conditions of $\mathbf{P}$ to $\bar{\mathbf{P}}$ as $\epsilon \to 0$ is proved in Section 3.

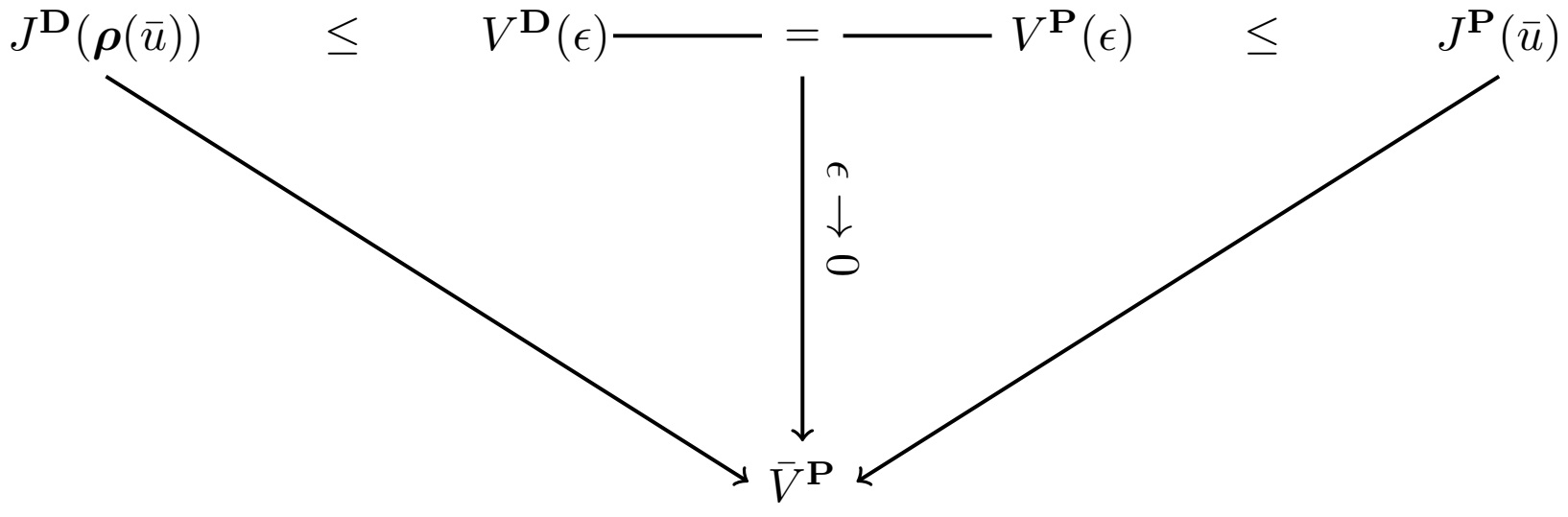

FIG. 1. *Methodology for obtaining the asymptotic error bound*

THEOREM 1. *The solution $V^{\mathbf{P}}(\epsilon)$ of the primal problem $\mathbf{P}$ and solution $V^{\mathbf{D}}(\epsilon)$ of the dual problem $\mathbf{D}$ satisfy the following equality*

$$V^{\mathbf{P}}(\epsilon) = V^{\mathbf{D}}(\epsilon), \quad \textit{for } \epsilon \in (0, \epsilon^*].$$

The strong duality property in Theorem 1 is paramount in obtaining the asymptotic error bound in equation (2) (see Fig.1). Furthermore, as the dual objective functional $J^{\mathbf{D}}$ is not necessarily coercive over its feasible set, the strong duality result, along with the existence of a unique solution to $\mathbf{P}$, will imply that a unique solution exists for $\mathbf{D}$.

THEOREM 2. *(a) Let $\bar{u}$ denote the optimal control for $\bar{\mathbf{P}}$. The control $\bar{u}$ provides an asymptotically optimal upper bound to the solution of $\mathbf{P}$ in the sense that*

$$V^{\mathbf{P}}(\epsilon) = J^{\mathbf{P}}(\bar{u}, \hat{z}, \epsilon) + O(\epsilon), \quad \textit{where } O(\epsilon) < 0 \ \textit{ as } \ \epsilon \to 0^+. \tag{11}$$

*In (11), $V^{\mathbf{P}}(\epsilon)$ is the solution to $\mathbf{P}$, $J^{\mathbf{P}}$ is the objective functional in (3), and $\hat{z}$ is the state satisfying the differential equations and boundary conditions in $\mathbf{P}$ with control given by $\bar{u}$.*

*(b) Consider the following control*

$$\boldsymbol{\rho}(t, \epsilon) = Q(t, \epsilon)\hat{z}(t, \epsilon), \tag{12}$$

*where $\hat{z}$ is the state satisfying the differential equations and boundary conditions in $\mathbf{P}$ with control given by $\bar{u}$. The control $\boldsymbol{\rho}$ provides an asymptotically optimal lower bound to the solution of $\mathbf{D}$ in the sense that*

$$V^{\mathbf{D}}(\epsilon) = J^{\mathbf{D}}(\boldsymbol{\rho}, \hat{\gamma}, \epsilon) + O(\epsilon), \quad \textit{where } O(\epsilon) > 0 \ \textit{ as } \ \epsilon \to 0^+. \tag{13}$$

*In (13), $V^{\mathbf{D}}(\epsilon)$ is the solution to $\mathbf{D}$, $J^{\mathbf{D}}$ is the objective functional in (6), and $\hat{\gamma}$ satisfies the differential equations in $\boldsymbol{D}$ with control given by $\boldsymbol{\rho}$ and boundary condition*

$$\hat{\gamma}(1, \epsilon) = -I^{\frac{1}{\epsilon}}\pi(\epsilon)Q(1, \epsilon)^{-1}\boldsymbol{\rho}(1, \epsilon). \tag{14}$$

*(c) From Theorem 1 and parts (a) and (b) in Theorem 2, we obtain the result*

$$J^{\mathbf{D}}(\boldsymbol{\rho}, \hat{\gamma}, \epsilon) \leq V^{\mathbf{D}}(\epsilon) = V^{\mathbf{P}}(\epsilon) \leq J^{\mathbf{P}}(\bar{u}, \hat{z}, \epsilon),$$

$$\left| J^{\mathbf{P}}(\bar{u}, \hat{z}, \epsilon) - J^{\mathbf{D}}(\boldsymbol{\rho}, \hat{\gamma}, \epsilon) \right| = O(\epsilon),$$

*as $\epsilon \to 0$ with $\hat{z}$ and $\hat{\gamma}$ as in Theorem 2, parts (a) and (b) respectively. Thus the constant $C$ in (2) is given by*

$$C = \frac{1}{\bar{V}^{\mathbf{P}}} \left( \lim_{\epsilon \to 0} \left| \frac{J^{\mathbf{P}}(\bar{u}, \hat{z}, \epsilon) - J^{\mathbf{D}}(\boldsymbol{\rho}, \hat{\gamma}, \epsilon)}{\epsilon} \right| \right). \tag{15}$$

Section 3 will be devoted to proving the convergence of the primal solution to the reduced solution (see Fig. 1). Theorems 1 and 2 will be proved in Section 4. Section 5 is reserved for showing the construction of the dual problem **D** and we give some numerical examples in Section 6.

**3. Convergence of the optimality conditions for the primal problem.** We begin this section with an outline of the necessary optimality conditions for **P** and its corresponding reduced problem $\bar{\mathbf{P}}$ before proceeding to show that the optimality conditions for **P** converge to those of $\bar{\mathbf{P}}$ as $\epsilon \to 0$. Omitting dependence on $t$ and $\epsilon$ for simplicity, the Hamiltonian function associated with an SPOC problem of the form **P** (see [18], Chap. 6) is defined as

$$H^{\mathbf{P}}(\hat{z}, \hat{u}, \hat{\chi}) = \frac{1}{2}(\hat{z}^T Q \hat{z} + \hat{u}^T R \hat{u}) + \hat{\chi}^T (A\hat{z} + I^{\epsilon} b \hat{u}),$$

where $\hat{\chi} \in W^{1,2}([0,1]; \mathbb{R}^{m+n})$ as a function of $t$ is the co-state variable. The scaling $\hat{\chi} \to I^{\frac{1}{\epsilon}} \hat{\chi}$ recovers the standard form of the Hamiltonian. Since there are no boundary conditions at the terminal time we have a normal Hamiltonian multiplier. Let $u$, $z$, $\chi$ denote the optimal control, state and co-state respectively of **P**. These variables must satisfy the following necessary optimality conditions (see [9])

$$\frac{dz}{dt} = \frac{dH^{\mathbf{P}}}{dI^{\epsilon}\chi}, \qquad \frac{dI^{\epsilon}\chi}{dt} = \frac{dH^{\mathbf{P}}}{dz}, \tag{16}$$

with boundary conditions

$$z(0, \epsilon) = z_0, \quad \chi(1, \epsilon) = I^{\frac{1}{\epsilon}} \pi(\epsilon) z(1, \epsilon). \tag{17}$$

Although the control $\hat{u}$ is constrained, we can obtain an explicit form for the optimal control, $u$, using Pontryagin's minimum principle. This principles states that $u$ must satisfy the following inequality

$$H^{\mathbf{P}}(z, u, \chi) \leq H^{\mathbf{P}}(z, \hat{u}, \chi), \qquad \text{for all} \quad \hat{u} \in U \tag{18}$$

From (18), it follows that

$$(u + f)^T R (u + f) \leq (\hat{u} + f)^T R (\hat{u} + f),$$

where $f = R^{-1} b^T I^{\epsilon} \chi$. Thus, for all $j = 1, \dots k$, the optimal control is given by

$$u_j = \begin{cases} -(R^{-1} b^T I^{\epsilon} \chi)_j & \text{if } \alpha_j \leq -(R^{-1} b^T I^{\epsilon} \chi)_j \leq \beta_j, \\ \alpha_j & \text{if } \alpha_j > -(R^{-1} b^T I^{\epsilon} \chi)_j, \\ \beta_j & \text{if } \beta_j < -(R^{-1} b^T I^{\epsilon} \chi)_j, \end{cases} \tag{19}$$

where the co-state $\chi$ must satisfy the relevant optimality conditions in (16), and (17).

Remark 5. *As the solution to* $\mathbf{P}$ *is unique, the necessary optimality conditions* (16) *and* (17) *are also sufficient and any solution satisfying these conditions will be the unique solution to* $\mathbf{P}$.

We now derive the optimality conditions for $\bar{\mathbf{P}}$. Omitting dependence on $t$, the Hamiltonian associated with $\bar{\mathbf{P}}$ is defined as

$$H^{\bar{\mathbf{P}}}(\hat{x}, \hat{u}, \hat{\chi}_1) = \frac{1}{2}(\hat{x}^T \mathcal{Q}\hat{x} + \hat{u}^T R^0 \hat{u}) + \hat{\chi}_1^T(\mathcal{A}\hat{x} + b_1^0 \hat{u}).$$

where $\hat{\chi}_1 \in W^{1,2}([0,1]; \mathbb{R}^m)$ as a function of $t$ is the co-state variable. Since there is no terminal condition on the problem $\bar{\mathbf{P}}$, the Hamiltonian has a normal multiplier. Let $\bar{u}$, $\bar{x}$ and $\bar{\chi}_1$ denote the optimal control, state and co-state respectively of $\bar{\mathbf{P}}$. These variables must satisfy the following necessary optimality conditions

$$\begin{aligned} \frac{d\bar{x}}{dt} &= \mathcal{A}(t)\bar{x}(t) + b_1^0(t)\bar{u}(t), \\ \frac{d\bar{\chi}_1}{dt} &= -\mathcal{A}^T(t)\bar{\chi}_1(t) - \mathcal{Q}(t)\bar{x}(t), \end{aligned} \tag{20}$$

with boundary conditions

$$\bar{x}(0) = z_{1,0}, \qquad \bar{\chi}_1(1) = \pi_{11}^0 \bar{x}(1). \tag{21}$$

Using Pontryagin's minimum principle, the optimal control is given by

$$\bar{u}_j(t) = \begin{cases} (-R^0(t)^{-1}b_1^0(t)^T\bar{\chi}_1(t))_j & \text{if } \alpha_j(t) \le (-R^0(t)^{-1}b_1^0(t)^T\bar{\chi}_1(t))_j \le \beta_j(t), \\ \alpha_j(t) & \text{if } \alpha_j(t) > (-R^0(t)^{-1}b_1^0(t)^T\bar{\chi}_1(t))_j, \\ \beta_j(t) & \text{if } \beta_j(t) < (-R^0(t)^{-1}b_1^0(t)^T\bar{\chi}_1(t))_j, \end{cases}$$

for all $j = 1, \ldots k$ and where $\bar{\chi}_1$ must satisfy the relevant conditions in (20) and (21).

Remark 6. *As the solution to* $\bar{\mathbf{P}}$ *is unique, the necessary optimality conditions are also sufficient and any solution satisfying these conditions will be the unique solution to* $\bar{\mathbf{P}}$.

Proposition 3. *The optimal states,* $z_1$ *and* $z_2$*, and respective co-states,* $\chi_1$ *and* $\chi_2$*, of* $\mathbf{P}$ *have the following asymptotic expansions*

$$\begin{aligned} z_1(t,\epsilon) &= z_{1,o}^0(t) + O(\epsilon), \\ \chi_1(t,\epsilon) &= \chi_{1,o}^0(t) + O(\epsilon), \\ z_2(t,\epsilon) &= z_{2,o}^0(t) + z_{2,i}^0(\tau) + O(\epsilon), \\ \chi_2(t,\epsilon) &= \chi_{2,o}^0(t) + \chi_{2,i}^0(\tau) + \chi_{2,f}^0(\sigma) + O(\epsilon), \end{aligned} \tag{22}$$

*as* $\epsilon \to 0$*, uniformly on* $[0,1]$ *where* $\tau$ *and* $\sigma$ *are defined as*

$$\tau = \frac{t}{\epsilon} \quad \text{and} \quad \sigma = \frac{1-t}{\epsilon}.$$

*The variables* $z_{1,o}^0$, $\chi_{1,o}^0$, $z_{2,o}^0$, *and* $\chi_{2,o}^0$ *satisfy the following system*

$$\begin{aligned}
\frac{dz_{1,o}^0}{dt} &= A_{11}^0(t)z_{1,o}^0 + A_{12}^0(t)z_{2,o}^0 + b_1^0(t)u_o^0, \\
\frac{d\chi_{1,o}^0}{dt} &= -A_{11}^0(t)^T\chi_{1,o}^0 - A_{21}^0(t)^T\chi_{2,o}^0 - Q_{11}^0(t)z_{1,o}^0 - Q_{12}^0(t)z_{2,o}^0, \\
0 &= A_{21}^0(t)z_{1,o}^0 + A_{22}^0(t)z_{2,o}^0, \\
0 &= -A_{12}^0(t)^T\chi_{1,o}^0 - A_{22}^0(t)^T\chi_{2,o}^0 - Q_{21}^0(t)z_{1,o}^0 - Q_{22}^0(t)z_{2,o}^0,
\end{aligned} \tag{23}$$

*with* $u_o^0$ *defined as*

(24)

$$(u_o^0(t))_j = \begin{cases} -(R^0(t)^{-1}b_1^0(t)\chi_{1,o}^0(t))_j & \text{if } \alpha_j(t) \le -(R^0(t)^{-1}b_1^0(t)\chi_{1,o}^0(t))_j \le \beta_j(t), \\ \alpha_j(t) & \text{if } \alpha_j(t) > -(R^0(t)^{-1}b_1^0(t)\chi_{1,o}^0(t))_j, \\ \beta_j(t) & \text{if } \beta_j(t) < -(R^0(t)^{-1}b_1^0(t)\chi_{1,o}^0(t))_j, \end{cases}$$

*for all* $j = 1, \dots k$. *The variables* $z_{1,o}^0$ *and* $\chi_{1,o}^0$ *satisfy the following boundary conditions*

$$z_{1,o}(0) = z_{1,0}, \qquad \chi_{1,o}(1) = \pi_{11}^0 z_{1,o}^0(1). \tag{25}$$

*The variables* $z_{2,i}^0$, $\chi_{2,i}^0$, $\chi_{2,f}^0$ *satisfy the differential equations*

$$\begin{aligned}
\frac{dz_{2,i}^0}{d\tau} &= A_{22}^0(0)z_{2,i}^0, \\
\frac{d\chi_{2,i}^0}{d\tau} &= -A_{22}^0(0)^T\chi_{2,i}^0 - Q_{22}^0(0)z_{2,i}^0, \\
\frac{d\chi_{2,f}^0}{d\sigma} &= A_{22}^0(1)^T\chi_{2,f}^0 + Q_{22}^0(1)z_{2,f}^0,
\end{aligned} \tag{26}$$

*along with the boundary conditions*

$$\begin{aligned}
z_{2,i}^0(0) &= z_{2,0} - z_{2,o}^0(0), \\
\chi_{2,i}^0(\infty) &= 0, \\
\chi_{2,f}^0(0) &= \pi_{22}^0 z_{2,o}^0(1) - \chi_{2,o}^0(1).
\end{aligned} \tag{27}$$

The proof of Proposition 3 will use the following lemma.

LEMMA 4. *Consider a fundamental solution matrix* $Z(t, \epsilon)$ *of the system*

$$\frac{d}{dt}\begin{bmatrix} \hat{R}_{z_1} \\ \hat{R}_{\chi_1} \\ \epsilon\hat{R}_{z_2} \\ \epsilon\hat{R}_{\chi_2} \end{bmatrix} = \begin{bmatrix} \mathbf{A}_{11}(t,\epsilon) & \mathbf{A}_{12}(t,\epsilon) \\ \mathbf{A}_{21}(t,\epsilon) & \mathbf{A}_{22}(t,\epsilon) \end{bmatrix} \begin{bmatrix} \hat{R}_{z_1} \\ \hat{R}_{\chi_1} \\ \hat{R}_{z_2} \\ \hat{R}_{\chi_2} \end{bmatrix}.$$

*where* $\mathbf{A}_{11}, \mathbf{A}_{12}, \mathbf{A}_{21}$ *and* $\mathbf{A}_{22}$ *are defined as*

(28)

$$\mathbf{A}_{11}(t,\epsilon) = \begin{bmatrix} A_{11}(t,\epsilon) & 0 \\ -Q_{11}(t,\epsilon) & -A_{11}(t,\epsilon)^T \end{bmatrix}, \quad \mathbf{A}_{12}(t,\epsilon) = \begin{bmatrix} A_{12}(t,\epsilon) & 0 \\ -Q_{12}(t,\epsilon) & -A_{21}(t,\epsilon)^T \end{bmatrix},$$

$$\mathbf{A}_{21}(t,\epsilon) = \begin{bmatrix} A_{21}(t,\epsilon) & 0 \\ -Q_{21}(t,\epsilon) & -A_{12}(t,\epsilon)^T \end{bmatrix}, \quad \mathbf{A}_{22}(t,\epsilon) = \begin{bmatrix} A_{22}(t,\epsilon) & 0 \\ -Q_{22}(t,\epsilon) & -A_{22}(t,\epsilon)^T \end{bmatrix}.$$

*The matrix $Z$ satisfies the bounds*

$$Z(t,\epsilon) = \left[\begin{array}{cc|cc} O(1) & O(1) & O(\epsilon e^{-\frac{kt}{\epsilon}}) & O(\epsilon e^{-\frac{k(1-t)}{\epsilon}}) \\ O(1) & O(1) & O(\epsilon e^{-\frac{kt}{\epsilon}}) & O(\epsilon e^{-\frac{k(1-t)}{\epsilon}}) \\ \hline O(1) & O(1) & O(e^{-\frac{kt}{\epsilon}}) & O(e^{-\frac{k(1-t)}{\epsilon}}) \\ O(1) & O(1) & O(e^{-\frac{kt}{\epsilon}}) & O(e^{-\frac{k(1-t)}{\epsilon}}) \end{array}\right], \tag{29}$$

(upper-left block $2m \times 2m$, upper-right block $2m \times 2n$, lower-left block $2n \times 2m$, lower-right block $2n \times 2n$)

*as $\epsilon \to 0$, uniformly on $[0,1]$ for some constant $k > 0$.*

*Proof.* Consider the following non-singular linear transformation (see [8])

$$\begin{bmatrix} \hat{R}_{z_1} \\ \hat{R}_{\chi_1} \\ \hat{R}_{z_2} \\ \hat{R}_{\chi_2} \end{bmatrix} = \begin{bmatrix} I_{2m} & -\epsilon H_1(t,\epsilon) \\ -H_2(t,\epsilon) & I_{2n} + \epsilon H_2(t,\epsilon) H_1(t,\epsilon) \end{bmatrix} \begin{bmatrix} \eta_{z_1} \\ \eta_{\chi_1} \\ \eta_{z_2} \\ \eta_{\chi_2} \end{bmatrix}. \tag{30}$$

The matrices $H_1$ and $H_2$ satisfy the equations

(31)

$$\epsilon \frac{dH_1}{dt} = \epsilon(\mathbf{A}_{11}(t,\epsilon) - \mathbf{A}_{12}(t,\epsilon)H_2)H_1 - H_1(\mathbf{A}_{22}(t,\epsilon) + \epsilon H_2 \mathbf{A}_{12}(t,\epsilon)) - \mathbf{A}_{12}(t,\epsilon),$$
$$\epsilon \frac{dH_2}{dt} = \mathbf{A}_{22}(t,\epsilon)H_2 - \epsilon H_2 \mathbf{A}_{11}(t,\epsilon) + \epsilon H_2 \mathbf{A}_{12}(t,\epsilon) H_2 - \mathbf{A}_{21}(t,\epsilon).$$

Upon specifying any arbitrary boundary conditions for $H_1$ and $H_2$ that are bounded as $\epsilon \to 0$, a well known result ([8] and [26] Lemma 9.2.1-9.2.2) guarantees that the solutions to (31) exist and are bounded as $\epsilon \to 0$, uniformly on $[0,1]$. Applying the transformation (30), the variables $\eta_{z_1}$, $\eta_{\chi_1}$, $\eta_{z_2}$ and $\eta_{\chi_2}$ satisfy the diagonalised system

$$\frac{d}{dt}\begin{bmatrix} \eta_{z_1} \\ \eta_{\chi_1} \end{bmatrix} = (\mathbf{A}_{11}(t,\epsilon) - \mathbf{A}_{12}(t,\epsilon)H_2(t,\epsilon)) \begin{bmatrix} \eta_{z_1} \\ \eta_{\chi_1} \end{bmatrix}, \tag{32}$$

$$\epsilon\frac{d}{dt}\begin{bmatrix} \eta_{z_2} \\ \eta_{\chi_2} \end{bmatrix} = (\mathbf{A}_{22}(t,\epsilon) + \epsilon H_2(t,\epsilon)\mathbf{A}_{12}(t,\epsilon)) \begin{bmatrix} \eta_{z_2} \\ \eta_{\chi_2} \end{bmatrix}. \tag{33}$$

Let $U_1(t,\epsilon)$ denote a fundamental solution matrix for the system (32). As all matrices in (32) are bounded as $\epsilon \to 0$ uniformly on $[0,1]$ and the system is unperturbed, it follows that $U_1$ satisfies the bound

$$U_1(t,\epsilon) = O(1), \tag{34}$$

as $\epsilon \to 0$, uniformly on $[0,1]$. To determine bounds on a fundamental solution matrix for (33), let us rewrite the system as

$$\epsilon\frac{d}{dt}\begin{bmatrix} \eta_{z_2} \\ \eta_{\chi_2} \end{bmatrix} = \Big(\mathbf{A}_{22}^0(t) + \epsilon T(t,\epsilon)\Big) \begin{bmatrix} \eta_{z_2} \\ \eta_{\chi_2} \end{bmatrix}, \tag{35}$$

where

$$T(t,\epsilon) = \left(\sum_{k=1}^{\infty} \epsilon^{k-1}\mathbf{A}_{22}^k(t) + H_2(t,\epsilon)\mathbf{A}_{12}(t,\epsilon)\right).$$

To find the solutions of (35), consider the transformation

(36)
$$\begin{bmatrix}\eta_{z_2}\\ \eta_{\chi_2}\end{bmatrix} = \begin{bmatrix} I_n & 0\\ \Omega(t) & I_n\end{bmatrix}\begin{bmatrix}\nu_{z_2}\\ \nu_{\chi_2}\end{bmatrix}$$

where $\Omega(t) \in W^{1,2}([0,1];\mathbb{R}^{n\times n})$ as a function of $t$ and satisfies the equation

(37)
$$\Omega(t)A_{22}^0(t) + Q_{22}^0(t) + A_{22}^0(t)^T\Omega(t) = 0,$$

for all $t \in [0,1]$. The equation (35) is transformed into the system

(38)
$$\epsilon\frac{d}{dt}\begin{bmatrix}\nu_{z_1}\\ \nu_{\chi_2}\end{bmatrix} = \left(\begin{bmatrix}A_{22}^0(t) & 0\\ 0 & -A_{22}^0(t)^T\end{bmatrix} + \epsilon S(t,\epsilon)\right)\begin{bmatrix}\nu_{z_1}\\ \nu_{\chi_2}\end{bmatrix},$$

where $S$ is uniformly bounded on $[0,1]$ as $\epsilon \to 0$. From a well known result in [14], the fundamental solution to a system of the form (38) satisfies the bounds

(39)
$$V(t,\epsilon) = \begin{bmatrix} O(e^{\frac{-kt}{\epsilon}}) & O(e^{\frac{-k(1-t)}{\epsilon}})\\ O(e^{\frac{-kt}{\epsilon}}) & O(e^{\frac{-k(1-t)}{\epsilon}})\end{bmatrix},$$

as $\epsilon \to 0$ uniformly on $[0,1]$. It follows from (30) and (36) that the fundamental solution $Z(t,\epsilon)$ can be written as

(40)
$$Z(t,\epsilon) = \begin{bmatrix} I & -\epsilon H_1(t,\epsilon)\\ -H_2(t,\epsilon) & I+\epsilon H_2(t,\epsilon)H_1(t,\epsilon)\end{bmatrix}\begin{bmatrix}U_1(t,\epsilon) & 0\\ 0 & U_2(t,\epsilon)\end{bmatrix},$$

where

$$U_2(t,\epsilon) = \begin{bmatrix} I_n & 0\\ \Omega(t) & I_n\end{bmatrix} V(t,\epsilon).$$

A multiplication of the matrices in (40) along with the bounds in (34) and (39) gives the bounds in (29). □

*Proof of Proposition 3.* Let us write the optimal variables of **P** as follows

(41)
$$\begin{aligned} z_1(t,\epsilon) &= z_{1,o}^0(t) + R_{z_1}(t,\epsilon), & z_2(t,\epsilon) &= z_{2,o}^0(t) + z_{2,i}^0(\tau) + R_{z_2}(t,\epsilon),\\ \chi_1(t,\epsilon) &= \chi_{1,o}^0(t) + R_{\chi_1}(t,\epsilon), & \chi_2(t,\epsilon) &= \chi_{2,o}^0(t) + \chi_{2,i}^0(\tau) + \chi_{2,f}^0(\sigma) + R_{\chi_2}(t,\epsilon),\end{aligned}$$

for all $\epsilon \in (0,\epsilon^*]$ and $t \in [0,1]$, where $R_{z_1}, R_{\chi_1}, R_{z_2}$, and $R_{\chi_2}$ are the remainders of the expansions. It follows from (26) and (27) along with assumption (a) that the variables $z_{2,i}^0$, $\chi_{2,i}^0$, and $\chi_{2,f}^0$ satisfy the following equalities

(42)
$$z_{2,i}^0(\tau) = O(e^{-k\tau}), \quad \chi_{2,i}^0(\tau) = O(e^{-k\tau}), \quad \chi_{2,f}^0(\sigma) = O(e^{-k\sigma}),$$

as $\tau,\sigma \to \infty$ and where $k$ is some fixed positive constant (see [14]). Furthermore, all outer layer terms are bounded on $[0,1]$. To prove Proposition 3 we must show

(43)
$$R_{z_1}(t,\epsilon), R_{\chi_1}(t,\epsilon), R_{z_2}(t,\epsilon), R_{\chi_2}(t,\epsilon) = O(\epsilon),$$

as $\epsilon \to 0$, uniformly on $[0,1]$. To do so, we first construct a system of differential equations that the remainder terms satisfy and then bound its solution. From equations (16), (23), (24), (26) and(41), we obtain the following system for the remainder

terms

$$\frac{d}{dt}\begin{bmatrix} R_{z_1} \\ R_{\chi_1} \\ \epsilon R_{z_2} \\ \epsilon R_{\chi_2} \end{bmatrix} = \begin{bmatrix} \mathbf{A}_{11}(t,\epsilon) & \mathbf{A}_{12}(t,\epsilon) \\ \mathbf{A}_{21}(t,\epsilon) & \mathbf{A}_{22}(t,\epsilon) \end{bmatrix} \begin{bmatrix} R_{z_1} \\ R_{\chi_1} \\ R_{z_2} \\ R_{\chi_2} \end{bmatrix} + \begin{bmatrix} \Gamma_{z_1}(t,\epsilon) \\ \Gamma_{\chi_1}(t,\epsilon) \\ \Gamma_{z_2}(t,\epsilon) \\ \Gamma_{\chi_2}(t,\epsilon) \end{bmatrix} + \begin{bmatrix} b_1(t,\epsilon)(u-u_o^0) \\ 0 \\ 0 \\ 0 \end{bmatrix}. \tag{44}$$

The matrices $\boldsymbol{A}_{ij}$, $i,j=1,2$ satisfy (28) and are bounded on $[0,1]$ as $\epsilon \to 0$. This, in addition to the bounds on the variables satisfying (23) and the bounds (42), yields the following bounds for $\Gamma_{z_1}, \Gamma_{\chi_1}, \Gamma_{z_2}$ and $\Gamma_{\chi_2}$

$$\begin{aligned} \Gamma_{z_1} &= O(\epsilon + e^{-k\tau}), \qquad & \Gamma_{\chi_1} &= O(\epsilon + e^{-k\tau} + e^{-k\sigma}), \\ \Gamma_{z_2} &= O(\epsilon + e^{-k\tau}), \qquad & \Gamma_{\chi_2} &= O(\epsilon + e^{-k\tau} + e^{-k\sigma}), \end{aligned} \tag{45}$$

as $\epsilon \to 0$ uniformly on $[0,1]$ for some constant $k > 0$. It follows from (41) evaluated at $t=0$ and $t=1$ along with the conditions, (17), (25) and (27) that the boundary conditions for (44) satisfy

$$L \begin{bmatrix} R_{z_1}(0,\epsilon) \\ R_{\chi_1}(0,\epsilon) \\ R_{z_2}(0,\epsilon) \\ R_{\chi_2}(0,\epsilon) \end{bmatrix} + R \begin{bmatrix} R_{z_1}(1,\epsilon) \\ R_{\chi_1}(1,\epsilon) \\ R_{z_2}(1,\epsilon) \\ R_{\chi_2}(1,\epsilon) \end{bmatrix} = B(\epsilon), \tag{46}$$

where

$$B(\epsilon)^T = \begin{bmatrix} 0_{1\times m} & O(\epsilon)_{1\times m} & 0_{1\times n} & O(\epsilon)_{1\times n} \end{bmatrix}, \tag{47}$$

as $\epsilon \to 0$ and $L$ and $R$ are block matrices of the form

$$L = \begin{bmatrix} I_m & 0 & 0 & 0 \\ 0 & 0 & 0 & 0 \\ 0 & 0 & I_n & 0 \\ 0 & 0 & 0 & 0 \end{bmatrix}, \qquad R = \begin{bmatrix} 0 & 0 & 0 & 0 \\ -\pi_{11}(\epsilon) & I_m & 0 & 0 \\ 0 & 0 & 0 & 0 \\ 0 & 0 & -\pi_{22}(\epsilon) & I_n \end{bmatrix}. \tag{48}$$

The solution to (44) and (46) is given by (see [26], Chap. 9.3)

$$\begin{bmatrix} R_{z_1}(t,\epsilon) \\ R_{\chi_1}(t,\epsilon) \\ R_{z_2}(t,\epsilon) \\ R_{\chi_2}(t,\epsilon) \end{bmatrix} = Z(t,\epsilon)M(\epsilon)^{-1}B(\epsilon) + \int_0^1 G(t,s,\epsilon) \begin{bmatrix} \Gamma_{z_1}(s) + b_1(s,\epsilon)(u(s,\epsilon) - u_o^0(s)) \\ \Gamma_{\chi_1}(s) \\ \frac{1}{\epsilon}\Gamma_{z_2}(s) \\ \frac{1}{\epsilon}\Gamma_{\chi_2}(s) \end{bmatrix} ds. \tag{49}$$

where $G$ is given by

$$G(t,s,\epsilon) = \begin{cases} Z(t,\epsilon)M(\epsilon)^{-1}L(\epsilon)Z(0,\epsilon)Z(s,\epsilon)^{-1} & \text{for } s < t, \\ -Z(t,\epsilon)M(\epsilon)^{-1}R(\epsilon)Z(1,\epsilon)Z(s,\epsilon)^{-1} & \text{for } s > t. \end{cases} \tag{50}$$

In (49) and (50), $M$ is defined as

$$M(\epsilon) = L(\epsilon)Z(0,\epsilon) + R(\epsilon)Z(1,\epsilon), \tag{51}$$

and $Z$ is a fundamental solution to the homogeneous system given in Lemma 4 satisfying the bounds in (29).

Remark 7. *Recall that we may specify arbitrary bounded as $\epsilon \to 0$ boundary conditions for $H_1$ satisfying (31). Hence, we choose boundary conditions that allow certain terms in the matrices $Z(0,\epsilon)$ and $Z(1,\epsilon)$ to be set to 0, thereby allowing us to obtain the $O(\epsilon)$ bounds for the remainder terms. We specify the following boundary conditions for $H_1$*

$$\begin{aligned}
&(H_1)_{11}(0,\epsilon) = 0_{m\times n}, && (H_1)_{12} = 0_{m\times n},\\
&(H_1)_{21}(1,\epsilon) = \pi_{11}(\epsilon)(H_1)_{11}(1,\epsilon), && (H_1)_{22} = \pi_{11}(\epsilon)(H_1)_{12}(1,\epsilon),
\end{aligned}$$

*With the above boundary conditions the terms $Z(0,\epsilon)$ and $Z(1,\epsilon)$ satisfy the bounds*

$$\begin{aligned}
Z(0,\epsilon) &= \left[\begin{array}{cc|cc} O(1) & O(1) & 0 & O(\epsilon e^{-\frac{k}{\epsilon}}) \\ O(1) & O(1) & O(\epsilon) & 0 \\ \hline O(1) & O(1) & O(1) & O(e^{-\frac{k}{\epsilon}}) \\ O(1) & O(1) & O(1) & O(e^{-\frac{k}{\epsilon}}) \end{array}\right], \\
Z(1,\epsilon) &= \left[\begin{array}{cc|cc} O(1) & O(1) & O(\epsilon e^{-\frac{k}{\epsilon}}) & 0 \\ O(1) & O(1) & 0 & O(\epsilon) \\ \hline O(1) & O(1) & O(e^{-\frac{k}{\epsilon}}) & O(1) \\ O(1) & O(1) & O(e^{-\frac{k}{\epsilon}}) & O(1) \end{array}\right],
\end{aligned} \tag{52}$$

*as $\epsilon \to 0$ for some constant $k > 0$.*

From equations (29), (47), (48), (51) and (52), it follows that the term $ZM^{-1}B$ satisfies the bound

$$Z(t,\epsilon)M(\epsilon)^{-1}B(\epsilon) = O(\epsilon), \tag{53}$$

as $\epsilon \to 0$ uniformly on $[0,1]$. Furthermore, the bounds (29) and (52) along with the definitions in (50) and (51) yield

$$G(t,s,\epsilon) = \left[\begin{array}{cc|cc} O(1) & O(1) & O(\epsilon) & O(\epsilon) \\ O(1) & O(1) & O(\epsilon) & O(\epsilon) \\ \hline O(1) & O(1) & O(\epsilon) & O(\epsilon) \\ O(1) & O(1) & O(\epsilon) & O(\epsilon) \end{array}\right], \tag{54}$$

as $\epsilon \to 0$, uniformly for $0 \le s < t \le 1$. Hence, from (45), (49), (53), and (54), it follows that

$$\begin{bmatrix} R_{z_1} \\ R_{\chi_1} \\ R_{z_2} \\ R_{\chi_2} \end{bmatrix} = \begin{bmatrix} O(\epsilon) \\ O(\epsilon) \\ O(\epsilon) \\ O(\epsilon) \end{bmatrix} + \int_0^1 G(t,s,\epsilon) \begin{bmatrix} b_1(s,\epsilon)(u(s,\epsilon) - u_o^0(s)) \\ 0 \\ 0 \\ 0 \end{bmatrix} ds. \tag{55}$$

From (19), (24) and (41), we obtain

$$||u(t,\epsilon) - u_o^0(t)||_\infty \le ||R^0(t)^{-1} b_1^0(t)^T R_{\chi_1}||_\infty + O(\epsilon), \quad \text{as } \epsilon \to 0, \tag{56}$$

uniformly on $[0,1]$ for all $j = 1,\dots,k$ where the norm is defined as

$$||x||_\infty = \max_{ij} |x_{ij}|, \tag{57}$$

for any matrix $x$. Substituting (56) into (55) and noting that $b_1$ is bounded on $[0,1]$ as $\epsilon \to 0$, an application of Gronwall's lemma yields

$$||R_{\chi_1}||_\infty = O(\epsilon) \quad \text{as } \epsilon \to 0,$$

uniformly on $[0,1]$. The bounds in (43) for the remainder terms $R_{z_1}$, $R_{z_2}$, and $R_{\chi_2}$ immediately follow upon integration. ∎

Corollary 5. *The optimal control has the asymptotic expansion*

$$u(t,\epsilon) = u_o^0(t) + O(\epsilon), \tag{58}$$

*as* $\epsilon \to 0$*, uniformly on* $[0,1]$.

*Proof.* The proof follows immediately from the equations (19), (22), and (24).

Proposition 6. *The optimal state* $z_1$*, co-state,* $\chi_1$ *and control* $u$ *of* $\mathbf{P}$ *converge uniformly on* $[0,1]$ *to the optimal state,* $\bar{x}$*, co-states* $\bar{\chi}_1$ *and control* $\bar{u}$ *respectively of* $\bar{\mathbf{P}}$ *as* $\epsilon \to 0$

*Proof.* The system (23) can be rewritten as

$$\begin{aligned} \frac{dz_1^0}{dt} &= \mathcal{A}(t) z_{1,o}^0 + b_1^0(t) u_o^0, \\ \frac{d\chi_{1,o}^0}{dt} &= -\mathcal{A}(t)^T \chi_{1,o}^0 - \mathcal{Q}(t) z_{1,o}^0, \end{aligned} \tag{59}$$

with $\mathcal{Q}$ and $\mathcal{A}$ defined in (9). From (24), (59), and the initial conditions (25), it is clear that $z_{1,o}^0$, $\chi_{1,o}^0$ and $u_o^0$ satisfy the optimality conditions of $\bar{\mathbf{P}}$. As the solution to $\bar{\mathbf{P}}$ is unique, it follows that

$$\begin{aligned} z_{1,o}^0(t) &= \bar{x}(t), \\ \chi_{1,o}^0(t) &= \bar{\chi}_1(t), \\ u_o^0(t) &= \bar{u}(t), \end{aligned} \tag{60}$$

on $[0,1]$. Substituting (60) into (22) and (58), we obtain the convergence result in Proposition 6. ∎

**4. Proofs of main theorems.** We begin this section by deriving the necessary optimality conditions for $\mathbf{D}$. Section 4.2 is dedicated to proving Theorem 1 and Section 4.2 to proving Theorem 2.

**4.1. Optimality conditions for the dual problem.** In a similar manner to Section 3, we begin by deriving the necessary optimality conditions for $\mathbf{D}$. Omitting dependence on $t$ and $\epsilon$ for simplicity, the Hamiltonian function associated with $\mathbf{D}$ is defined as

$$H^{\mathbf{D}}(\hat{\gamma}, \hat{\rho}, \hat{\mu}) = -\frac{1}{2}\hat{\rho}^T Q^{-1} \hat{\rho} - \theta(\hat{\gamma}) + \hat{\mu}^T(-A^T \hat{\gamma} + \hat{\rho}),$$

where $\hat{\mu} \in W^{1,2}([0,1]; \mathbb{R}^{m+n})$ as a function of $t$ is the co-state variable. The scaling $\hat{\mu} \to I^{\frac{1}{\epsilon}} \hat{\mu}$ recovers the standard form of the Hamiltonian. Since there are no boundary conditions at the initial and final time, we have a normal multiplier. Let $\rho$, $\gamma$, $\mu$ denote the optimal control, state and co-state respectively of **D**. The necessary optimality conditions which these variables must satisfy are given by

$$\frac{dI^\epsilon \mu}{dt} = A(t,\epsilon)\mu(t,\epsilon) + I^\epsilon b\xi(t,\epsilon), \qquad (61)$$
$$\frac{dI^\epsilon \gamma}{dt} = -A(t,\epsilon)^T\gamma(t,\epsilon) + \rho(t,\epsilon).$$

For all $j = 1, \ldots k$, $\xi$ is defined, component-wise, by

$$\xi_j = \begin{cases} (R^{-1}b^T I^\epsilon \gamma)_j & \text{if } \alpha_j \le (R^{-1}b^T I^\epsilon\gamma)_j \le \beta_j, \\ \alpha_j & \text{if } \alpha_j > (R^{-1}b^T I^\epsilon\gamma)_j, \\ \beta_j & \text{if } \beta_j < (R^{-1}b^T I^\epsilon\gamma)_j. \end{cases}$$

where we have forgone the dependence on $t$ and $\epsilon$ for simplicity. Using the calculus of variations, the optimal variables must satisfy the boundary conditions

$$\mu(0,\epsilon) = z_0 \qquad \mu(1,\epsilon) = -\pi(\epsilon)^{-1}I^\epsilon\gamma(1,\epsilon). \qquad (62)$$

As the optimal control in **D** is unconstrained, the Pontryagin maximum principle states that $\rho$ must satisfy the equality, $\frac{dH^{\mathbf{D}}}{d\rho} = 0$. Hence

$$\rho = Q(t,\epsilon)\mu,$$

where $\mu$ must satisfy the relevant optimality conditions given in (61) and (62).

**4.2. Proof of Theorem 1.** A strong duality property is needed for the primal and dual problem. Clearly this property is necessary for the asymptotic error bound (see Fig 1). More subtly, however, the property is needed to prove the uniqueness of the solution to **D**.

*Proof.* Suppose that $z$, $\chi$, and $u$ denote the optimal state, co-state and control respectively of **P**. Consider the following definitions for $\gamma$, $\mu$, and $\rho$

$$\gamma(t,\epsilon) = -\chi(t,\epsilon), \qquad (63)$$
$$\mu(t,\epsilon) = z(t,\epsilon),$$
$$\rho(t,\epsilon) = Q(t,\epsilon)z(t,\epsilon),$$

for $t \in [0,1]$ , $\epsilon \in (0,\epsilon^*]$. We first show that $(\gamma,\rho)$ is a feasible solution for the dual problem and then show that this solution is optimal and that strong duality holds. Substituting the definitions in (63) into the differential equation for $\chi$ in (16) and the boundary conditions (17) yields

$$I^\epsilon\dot{\gamma} = -A^T\gamma + \rho, \qquad (64)$$
$$\mu(0,\epsilon) = z_0,$$
$$\mu(1,\epsilon) = -\pi(\epsilon)^{-1}I^\epsilon\gamma(1,\epsilon).$$

The equations in (64) are equivalent to the equations for $\gamma$ in (61) with boundary conditions in (62) for **D**. Hence $(\gamma,\rho)$ is a feasible solution of the dual.

From weak duality, we know $V^{\mathbf{D}}(\epsilon) \le V^{\mathbf{P}}(\epsilon)$ (see [11], Chap. 2). To show that there is a zero duality gap, we need to show $V^{\mathbf{P}}(\epsilon) = V^{\mathbf{D}}(\epsilon)$. Evaluating (6) with the state and control given in (63) along with the substitution $\gamma(1) = -I^{\frac{1}{\epsilon}}\pi\mu(1) = -I^{\frac{1}{\epsilon}}\pi z(1)$ gives

$$
\begin{aligned}
&J^{\mathbf{D}}(Qz, -\chi, -\chi(0), -I^{\frac{1}{\epsilon}}\pi z(1), \epsilon)\\
&\quad= \int_0^1 \Big( -\frac{1}{2}z^TQz - \theta(-\chi)\Big)dt + \chi(0)^TI^{\epsilon}z_0 - \frac{1}{2}z(1)^T\pi z(1),\\
&\quad= \int_0^1 \Big( -\frac{1}{2}z^TQz - \theta(-\chi)\Big)dt - \int_0^1 \frac{d}{dt}\langle I^{\epsilon}\chi, z\rangle dt + \frac{1}{2}z(1)^T\pi z(1),\\
&\quad= \int_0^1 \Big( -\frac{1}{2}z^TQz - \theta(-\chi) - \langle I^{\epsilon}\dot{\chi}, z\rangle - \langle \chi, I^{\epsilon}\dot{z}\rangle\Big)dt + \frac{1}{2}z(1)^T\pi z(1).
\end{aligned}
\tag{65}
$$

From the differential equations in (16), we can evaluate the inner products in (65),

$$
\begin{aligned}
\langle I^{\epsilon}\dot{\chi}, z\rangle &= \langle -A^T\chi - Qz, z\rangle,\\
&= -\langle A^T\chi, z\rangle - \langle Qz, z\rangle,\\
\langle \chi, I^{\epsilon}\dot{z}\rangle &= \langle \chi, Az + I^{\epsilon}bu\rangle,\\
&= \langle A^T\chi, z\rangle + \langle b^TI^{\epsilon}\chi, u\rangle.
\end{aligned}
\tag{66}
$$

Substituting (66) into (65), gives

$$
J^{\mathbf{D}}(\epsilon) = \int_0^1 \Big(\frac{1}{2}z^TQz - \theta(-\chi) - \langle b^TI^{\epsilon}\chi, u\rangle\Big)dt + \frac{1}{2}z(1)^T\pi z(1),
\tag{67}
$$

Since $u$ is the optimal control given by (19), the inner product becomes

$$
\langle b^TI^{\epsilon}\chi, u\rangle_j = \begin{cases} -\frac{1}{R_{jj}}(\chi^TI^{\epsilon}b)_j^2 & \text{if } \alpha_j \le -(R^{-1}b^TI^{\epsilon}\chi)_j \le \beta_j,\\ (\chi^TI^{\epsilon}b)_j\alpha_j & \text{if } \alpha_j > -(R^{-1}b^TI^{\epsilon}\chi)_j,\\ (\chi^TI^{\epsilon}b)_j\beta_j & \text{if } \beta_j < -(R^{-1}b^TI^{\epsilon}\chi)_j, \end{cases}
$$

for all $j = 1, \dots k$. Hence $-\theta(-\chi) - \langle b^TI^{\epsilon}\chi, u\rangle = \sum_j^k \Omega_j(\chi)$ where $\Omega_j(\chi)$ is given by

$$
\Omega_j(\chi) = \begin{cases} \frac{1}{2R_{jj}}(\chi^TI^{\epsilon}b)_j^2 & \text{if } \alpha_j \le -(R^{-1}b^TI^{\epsilon}\chi)_j \le \beta_j,\\ \frac{1}{2}\alpha_j^2R_{jj} & \text{if } \alpha > -(R^{-1}b^TI^{\epsilon}\chi)_j,\\ \frac{1}{2}\beta_j^2R_{jj} & \text{if } \beta < -(R^{-1}b^TI^{\epsilon}\chi)_j, \end{cases}
\tag{68}
$$

for all $j = 1, \dots k$. Substituting (68) into (67), we obtain

$$
J^{\mathbf{D}}(Qz, -I^{\frac{1}{\epsilon}}\pi z(1), \epsilon) = \int_0^1 \Big(\frac{1}{2}z^TQz + \frac{1}{2}uRu\Big)dt + \frac{1}{2}z(1)^T\pi z(1) = V^{\mathbf{P}}(\epsilon).
$$

Since $(\gamma, \rho)$ is a feasible solution, by weak duality, we must have that $(\gamma, \rho)$ is optimal and $V^{\mathbf{D}}(\epsilon) = V^{\mathbf{P}}(\epsilon)$ ☐

As the solution to **P** is unique, we can justify Remark 3, i.e. that a unique solution exists for **D**.

**4.3. Proof of Theorem 2.** We split the proof of Theorem 2 into two subsections for parts (a) and (b) respectively. Part (c) follows immediately from parts (a) and (b) along with Theorem 1.

**4.3.1. Proof of Theorem 2 part (a).** The proof of Theorem 2 part (a) follows from a simple integration of the differential equations in $\mathbf{P}$ with the optimal control of the reduced problem $\bar{u}$. From the definition given in (3), we obtain

$$\begin{aligned}(69)\quad & |V^{\mathbf{P}}(\epsilon) - J^{\mathbf{P}}(\bar{u}, \hat{z}, \epsilon)| \\ & = \left| \frac{1}{2}\int_0^1 z(t,\epsilon)^T Q(t,\epsilon) z(t,\epsilon) + u(t,\epsilon)^T R(t,\epsilon) u(t,\epsilon) dt + \frac{1}{2} z(1,\epsilon)^T \pi(\epsilon) z(1,\epsilon) \right. \\ & \left. -\frac{1}{2}\int_0^1 \hat{z}(t,\epsilon)^T Q(t,\epsilon)\hat{z}(t,\epsilon) + \bar{u}(t)^T R(t,\epsilon)\bar{u}(t) dt - \frac{1}{2}\hat{z}(1,\epsilon)^T \pi(\epsilon)\hat{z}(1,\epsilon) \right|,\end{aligned}$$

where $\hat{z}$ solves the differential equations in $\mathbf{P}$ with control given by $\bar{u}$. Using the variation of parameters technique, we may write $\hat{z}$ and $z$ respectively as

$$(70)\quad \begin{aligned}\hat{z}(t,\epsilon) &= \Phi_{I\frac{1}{\epsilon}A}(t,0,\epsilon) z_0 + \int_0^t \Phi_{I\frac{1}{\epsilon}A}(t,s,\epsilon) b(s,\epsilon)\bar{u}(s) ds, \\ z(t,\epsilon) &= \Phi_{I\frac{1}{\epsilon}A}(t,0,\epsilon) z_0 + \int_0^t \Phi_{I\frac{1}{\epsilon}A}(t,s,\epsilon) b(s,\epsilon) u(s,\epsilon) ds,\end{aligned}$$

where $\Phi_{I\frac{1}{\epsilon}A}$ is the resolvent matrix for the differential equations in $\mathbf{P}$. Note that for a resolvent matrix $\Phi_x$ with $x \in \mathbb{R}^{(m+n)\times(m+n)}$, the following conditions must be satisfied for all $t \in [0,1]$

1. $\frac{d\Phi_x}{dt} = x(t,\epsilon)\Phi_x$,
2. $\Phi_x(t,t,\epsilon) = I_{m+n}$,
3. $\det(\Phi_x) \neq 0$.

Let us partition the matrix $\Phi_{I\frac{1}{\epsilon}A}$ as follows

$$\Phi_{I\frac{1}{\epsilon}A}(t,s,\epsilon) = \begin{bmatrix}\phi_{11}(t,s,\epsilon) & \phi_{12}(t,s,\epsilon) \\ \phi_{21}(t,s,\epsilon) & \phi_{22}(t,s,\epsilon)\end{bmatrix}, \quad \phi_{11} \in \mathbb{R}^{m\times m}, \ \ \phi_{22} \in \mathbb{R}^{n\times n},$$

A well known result (see [14] and [26], Theorem 6.1.2) guarantees that the matrices $\phi_{ij}$ for $i,j = 1,2$ satisfy the following bounds

$$(71)\quad \begin{aligned} ||\phi_{11}(t,s,\epsilon)||_\infty &\le M_{11}, \quad ||\phi_{12}(t,s,\epsilon)||_\infty \le \epsilon M_{12}, \\ ||\phi_{21}(t,s,\epsilon)||_\infty &\le M_{21}, \quad ||\phi_{22}(t,s,\epsilon)||_\infty \le (\epsilon + e^{\frac{-k(t-s)}{\epsilon}}) M_{22},\end{aligned}$$

uniformly for all $0 \le s \le t \le 1$, $\epsilon \in (0,\epsilon^*]$ where $k$ and $M_{ij}$ for $i,j = 1,2$ are fixed, positive constants and the norm is given in (57). As $b$ is bounded for all $0 \le t \le 1$ and $\epsilon \in (0,\epsilon^*]$, it follows from (58), (60), (70) and (71) that

$$(72)\quad z(t,\epsilon) = \hat{z}(t,\epsilon) + O(\epsilon) \quad \text{as} \ \ \epsilon \to 0.$$

Substituting (58), (60) and (72) into (69), gives the inequality in (11). As $u$ is the optimal control, it follows that

$$J^{\mathbf{P}}(u,z,\epsilon) \le J^{\mathbf{P}}(\bar{u},\hat{z},\epsilon).$$

Hence, we may conclude that that the $O(\epsilon)$ term in (11) satisfies $O(\epsilon) < 0$ as $\epsilon \to 0$.

**4.3.2. Proof of Theorem 2 part (b).** From the definition in (6), we obtain

(73)
$$\left|V^{\mathbf{D}} - J^{\mathbf{D}}(\boldsymbol{\rho},\hat{\gamma})\right| = \left| \int_0^1 -\frac{1}{2}\rho^T Q^{-1}\rho - \theta(\gamma)dt - \gamma(0)^T I^{\epsilon} z_0 - \frac{1}{2}\gamma(1)^T I^{\epsilon}\pi^{-1}I^{\epsilon}\gamma(1) \right.$$
$$\left. + \int_0^1 \frac{1}{2}\boldsymbol{\rho}^T Q^{-1}\boldsymbol{\rho} + \theta(\hat{\gamma})dt + \hat{\gamma}(0)^T I^{\epsilon} z_0 + \frac{1}{2}\hat{\gamma}(1)^T I^{\epsilon}\pi^{-1}I^{\epsilon}\hat{\gamma}(1)\right|,$$

where $\tilde{\gamma}$ is the state that satisfies the equations in $\boldsymbol{D}$ with control given by $\boldsymbol{\rho}$ in (12) and boundary condition (14). We omit the dependence of (73) on $\epsilon$ for simplicity. Let $\rho$ denote the optimal control of $\mathbf{D}$. It follows from (12), (63) and (72) that

$$\rho(t,\epsilon) = \boldsymbol{\rho}(t,\epsilon) + O(\epsilon), \qquad \text{as } \epsilon \to 0, \tag{74}$$

uniformly on $[0,1]$. Furthermore, from (14), (62), (63), and (74), it follows that

$$\gamma(1,\epsilon) = \hat{\gamma}(1,\epsilon) + O(\epsilon), \qquad \text{as } \epsilon \to 0, \tag{75}$$

uniformly on $[0,1]$. Let us define $\hat{\gamma}_\omega(\omega) = \hat{\gamma}(t)$ and $\gamma_\omega(\omega) = \gamma(t)$, where

$$\omega = 1 - t. \tag{76}$$

Using the variation of parameters technique and the we may write $\hat{\gamma}_\omega$ and $\gamma_\omega$ respectively as

$$\begin{aligned} \hat{\gamma}_\omega(\omega,\epsilon) &= \Phi_{I^{\frac{1}{\epsilon}}A^T}(\omega,0,\epsilon)\hat{\gamma}_\omega(0,\epsilon) - \int_0^\omega \Phi_{I^{\frac{1}{\epsilon}}A^T}(\omega,r,\epsilon) I^{\frac{1}{\epsilon}}\boldsymbol{\rho}(r,\epsilon)dr, \\ \gamma_\omega(\omega,\epsilon) &= \Phi_{I^{\frac{1}{\epsilon}}A^T}(\omega,0,\epsilon)\gamma_\omega(0,\epsilon) - \int_0^\omega \Phi_{I^{\frac{1}{\epsilon}}A^T}(\omega,r,\epsilon) I^{\frac{1}{\epsilon}}\rho(r,\epsilon)dr, \end{aligned} \tag{77}$$

where $\Phi_{I^{\frac{1}{\epsilon}}A^T}$ is the resolvent matrix for the differential equations in (61) under the transformation (76). The matrix $\Phi_{I^{\frac{1}{\epsilon}}A^T}$ satisfies the bounds in (71) for different fixed positive constants $k$ and $M_{ij}$ for $i,j = 1,2$ (see [14] and [26], Theorem 6.1.2). As $\Phi_{I^{\frac{1}{\epsilon}}A^T}(t,s,\epsilon)$ satisfies the bounds in (71) for $s < t$, it follows from (74) - (75) that

$$\gamma(t,\epsilon) = \hat{\gamma}(t,\epsilon) + O(\epsilon), \qquad \text{as } \epsilon \to 0, \tag{78}$$

uniformly on $[0,1]$.

Substituting (74) and (78) into (73), we obtain the inequality in (13). As $\rho$ is the optimal control, it follows that

$$J^{\mathbf{D}}(\rho,\gamma,\epsilon) \geq J^{\mathbf{D}}(\boldsymbol{\rho},\hat{\gamma},\epsilon).$$

Hence, we may conclude that that the $O(\epsilon)$ term in (11) satisfies $O(\epsilon) > 0$ as $\epsilon \to 0$.

**5. Dual construction.** In this section, we outline the steps used to construct the dual problem $\mathbf{D}$. Following the methods of [1], [6] and [7], we begin by converting the feasible set $\Sigma$ defined in (5) into a subspace. We introduce dummy variables $\hat{s}_0, \hat{s}_1 \in \mathbb{R}^{m+n}$ so that the problem $\mathbf{P}$ becomes

$$\begin{cases} \underset{u\in U}{\text{minimise}} & J^{\mathbf{P}}(\hat{z},\hat{u},\epsilon,\hat{s}_0,\hat{s}_1), \\ \text{subject to} & (\hat{z},\hat{u},\hat{s}_0,\hat{s}_1) \in \boldsymbol{\Sigma}, \end{cases}$$

where

$$J^{\mathbf{P}}(\hat{z},\hat{u},\epsilon,\hat{s}_0,\hat{s}_1)=\frac{1}{2}\int_0^1\Big(\hat{z}(t,\epsilon)^TQ(t,\epsilon)\hat{z}(t,\epsilon)+\hat{u}(t,\epsilon)^TR(t,\epsilon)\hat{u}(t,\epsilon)\Big)dt + \frac{1}{2}\hat{s}_1^T\pi(\epsilon)\hat{s}_1+\delta_{z_0}(\hat{s}_0)+\delta_U(\hat{u}),$$

and

$$\mathbf{\Sigma}=\Big\{(\hat{z},\hat{u},\hat{s}_0,\hat{s}_1):I^\epsilon\frac{d\hat{z}}{dt}=A(t,\epsilon)\hat{z}(t,\epsilon)+I^\epsilon b(t,\epsilon)\hat{u}(t,\epsilon),\hat{z}(0,\epsilon)=\hat{s}_0, \hat{z}(1,\epsilon)=\hat{s}_1,t\in[0,1],\epsilon\in(0,\epsilon^*]\Big\}.$$

The $\delta$-function is defined by

$$\delta_C(x)=\begin{cases}0 & \text{if } x\in C,\\ +\infty & \text{otherwise.}\end{cases}$$

The feasible set $\mathbf{\Sigma}$ is now a closed subspace of $W^{1,2}\times L^2\times\mathbb{R}^{n+m}\times\mathbb{R}^{n+m}$. Following standard techniques in duality theory, we introduce further dummy variables $\hat{v}\in W^{1,2}([0,1];\mathbb{R}^m,\mathbb{R}^n)$ , $\hat{w}\in L^2([0,1];\mathbb{R}^k)$, as functions of $t$, and $\hat{\kappa}_0,\hat{\kappa}_1\in\mathbb{R}^{n+m}$ in order to dualise the problem. The objective functional and conditions become

(79)
$$J^{\mathbf{P}}(\hat{v},\hat{w},\epsilon,\hat{\kappa}_0,\hat{\kappa}_1)=\frac{1}{2}\int_0^1\Big(\hat{v}(t,\epsilon)^TQ(t,\epsilon)\hat{v}(t,\epsilon)+\hat{w}(t,\epsilon)^TR(t,\epsilon)\hat{w}(t,\epsilon)\Big)dt + \frac{1}{2}\hat{\kappa}_1^T\pi(\epsilon)\hat{\kappa}_1+\delta_{z_0}(\hat{\kappa}_0)+\delta_U(\hat{w}),$$

subject to

$$\hat{v}(t,\epsilon)=\hat{z}(t,\epsilon),\quad \hat{w}(t,\epsilon)=\hat{u}(t,\epsilon),\quad \hat{\kappa}_0=\hat{s}_0,\quad \hat{\kappa}_1=\hat{s}_1,\quad (\hat{z},\hat{u},\hat{s}_0,\hat{s}_1)\in\mathbf{\Sigma}.$$

We separate the terms in (79) into the following four functions

$$\begin{aligned}f_1(\hat{v})&=\frac{1}{2}\int_0^1\hat{v}(t,\epsilon)^TQ(t,\epsilon)\hat{v}(t,\epsilon)dt,\\ f_2(\hat{w})&=\frac{1}{2}\int_0^1\hat{w}(t,\epsilon)^TR(t,\epsilon)\hat{w}(t,\epsilon)dt+\delta_U(\hat{w}),\\ f_3(\hat{\kappa}_0)&=\delta_{z_0}(\hat{\kappa}_0),\\ f_4(\hat{\kappa}_1)&=\frac{1}{2}\hat{\kappa}_1^T\pi(\epsilon)\hat{\kappa}_1.\end{aligned}$$

The dual functional $J^{\mathbf{D}}$ can be written by Fenchel duality [13] as

$$J^{\mathbf{D}}(\hat{\rho},\hat{\lambda}_2,\hat{\lambda}_3,\hat{\lambda}_4)=\begin{cases}-f_1^*(\hat{\rho})-f_2^*(\hat{\lambda}_2)-f_3^*(\hat{\lambda}_3)-f_4^*(\hat{\lambda}_4) & \text{if } (\hat{\rho},\hat{\lambda}_2,\hat{\lambda}_3,\hat{\lambda}_4)\in\mathbf{\Sigma}_1,\\ -\infty & \text{otherwise,}\end{cases}$$

where $\mathbf{\Sigma}_1$ is orthogonal to $\mathbf{\Sigma}$ and $f_1^*(\hat{\rho})$,$f_2^*(\hat{\lambda}_2)$, $f_3^*(\hat{\lambda}_3)$, and $f_4^*(\hat{\lambda}_4)$ are the Fenchel duals (see [6], [7]) of $f_1(\hat{v})$, $f_2(\hat{w})$,$f_3(\hat{\kappa}_0)$ and $f_4(\hat{\kappa}_1)$ respectively. The Fenchel duals

are defined as follows

$$
\begin{aligned}
f_1^*(\hat{\rho}) &= \sup_{\hat{v}(t,\epsilon)\in W^{1,2}} \int_0^1 \hat{\rho}(t,\epsilon)^T\hat{v}(t,\epsilon) - \frac{1}{2}\hat{v}(t,\epsilon)^T Q(t,\epsilon)\hat{v}(t,\epsilon)dt,\\
f_2^*(\hat{\lambda}_2) &= \sup_{\hat{w}(t,\epsilon)\in U} \int_0^1 \hat{\lambda}_2(t,\epsilon)^T\hat{w}(t,\epsilon) - \frac{1}{2}\hat{w}(t,\epsilon)^T R(t,\epsilon)\hat{w}(t,\epsilon)dt,\\
f_3^*(\hat{\lambda}_3) &= \sup_{\hat{\kappa}_0\in\mathbb{R}^{m+n}} \hat{\lambda}_3(\epsilon)^T\hat{\kappa}_0 - \delta_{z_0}(\hat{\kappa}_0),\\
f_4^*(\hat{\lambda}_4) &= \sup_{\hat{\kappa}_1\in\mathbb{R}^{m+n}} \hat{\lambda}_4(\epsilon)^T\hat{\kappa}_1 - \frac{1}{2}\hat{\kappa}_1^T\pi(\epsilon)\hat{\kappa}_1.
\end{aligned}
\tag{80}
$$

Hence, the dual problem can be written as

$$
\begin{cases}
\underset{\hat{\rho},\hat{\lambda}_2,\hat{\lambda}_3,\hat{\lambda}_4}{\text{maximise}} & J^{\mathbf{D}}(\hat{\rho},\hat{\lambda}_2,\hat{\lambda}_3,\hat{\lambda}_4,\epsilon),\\
\text{subject to} & (\hat{\rho},\hat{\lambda}_2,\hat{\lambda}_3,\hat{\lambda}_4)\in\mathbf{\Sigma}_1,
\end{cases}
\tag{81}
$$

where

$$
J^{\mathbf{D}} = -f_1^*(\hat{\rho}) - f_2^*(\hat{\lambda}_2) - f_3^*(\hat{\lambda}_3) - f_4^*(\hat{\lambda}_4).
$$

We will simplify the problem in (81) by finding the solutions to the functions in (80). We may rewrite $-f_1^*(\hat{\rho})$ as

$$
-f_1^*(\hat{\rho}) = \inf_{\hat{v}(t,\epsilon)\in W^{1,2}} \int_0^1 F(\hat{v},t,\epsilon)dt,
$$

where

$$
F(\hat{v},t,\epsilon) = \frac{1}{2}\hat{v}(t,\epsilon)^T Q(t,\epsilon)\hat{v}(t,\epsilon) - \hat{\rho}(t,\epsilon)^T\hat{v}(t,\epsilon).
\tag{82}
$$

Let $v(t,\epsilon)$ denote the optimal solution to (82). From calculus of variations, we know that $v(t,\epsilon)$ solves (82) if and only if it satisfies the Euler-Lagrange equation

$$
\frac{d}{dt}\frac{\partial F}{\partial \dot{v}} = \frac{\partial F}{\partial v}.
$$

Since $\frac{\partial F}{\partial \dot{v}} = 0$, we have $\frac{\partial F}{\partial v} = 0$ which implies $\rho(t,\epsilon) = Q(t,\epsilon)v(t,\epsilon)$. Therefore,

$$
f_1^*(\hat{\rho}) = \frac{1}{2}\int_0^1 \hat{\rho}(t,\epsilon)^T Q(t,\epsilon)^{-1}\hat{\rho}(t,\epsilon)dt.
\tag{83}
$$

For $-f_2^*(\hat{\lambda}_2)$ we rewrite the problem as

$$
\begin{cases}
\inf & l(1),\\
\text{subject to} & \dot{l} = \frac{1}{2}\hat{w}(t,\epsilon)^T R(t,\epsilon)\hat{w}(t,\epsilon) - \lambda_2(t,\epsilon)^T\hat{w}(t,\epsilon),\\
& \hat{w}(t,\epsilon)\in U, \quad l(0)=0,
\end{cases}
$$

where

$$
l(1) = \int_0^1 \frac{1}{2}\hat{w}(t,\epsilon)^T R(t,\epsilon)\hat{w}(t,\epsilon) - \hat{\lambda}_2(t,\epsilon)^T\hat{w}(t,\epsilon)dt.
$$

The Hamiltonian associated with $-f_2^*(\hat{\lambda}_2)$ is given by

$$\hat{H}(\hat{w},\hat{\mu},t,\epsilon) = \hat{\Delta}(t,\epsilon)^T(\frac{1}{2}\hat{w}(t,\epsilon)^T R(t,\epsilon)\hat{w}(t,\epsilon) - \hat{\lambda}_2(t,\epsilon)^T\hat{w}(t,\epsilon)).$$

Let $w$ and $\Delta$ denote the optimal control and co-state respectively. The variable $\Delta$ must satisfy the necessary optimality condition $\frac{d\Delta}{dt} = -\frac{d\hat{H}}{dl}$. As $\frac{d\hat{H}}{dl} = 0$, it follows that $\Delta(t,\epsilon) = c$ for some constant $c$ for all $t \in [0,1]$. By the transversality condition (see [7]), $\Delta(1) = 1$, hence

$$\Delta(t) = 1, \qquad \text{for all } t \in [0,1].$$

Applying Pontryagin's minimum principle, the optimal control $w$ is given by

$$w_j(t,\epsilon) = \begin{cases} (R(t,\epsilon)^{-1}\hat{\lambda}_2(t,\epsilon))_j & \text{if } \alpha_j(t) \le (R(t,\epsilon)^{-1}\hat{\lambda}_2(t,\epsilon))_j \le \beta_j(t), \\ \alpha_j(t) & \text{if } \alpha_j(t) > (R(t,\epsilon)^{-1}\hat{\lambda}_2(t,\epsilon))_j, \\ \beta_j(t) & \text{if } \beta_j(t) < (R(t,\epsilon)^{-1}\hat{\lambda}_2(t,\epsilon))_j, \end{cases} . \tag{84}$$

for all $j = 1, \ldots k$. Substituting (84) into $f_2^*(\hat{\lambda}_2)$, we obtain

$$f_2^*(\hat{\lambda}_2) = \int_0^1 \theta(\hat{\lambda}_2)dt, \tag{85}$$

where $\theta(\hat{\lambda}_2) = \sum_{j=1}^k \theta_j(\hat{\lambda}_2)$ and $\theta_j$ is defined by

$$\theta_j(\hat{\lambda}_2) = \begin{cases} \frac{1}{2R_{jj}(t,\epsilon)}(\hat{\lambda}_2(t,\epsilon))_j^2 & \text{if } \alpha_j(t) \le (R(t,\epsilon)^{-1}\hat{\lambda}_2(t,\epsilon))_j \le \beta_j(t), \\ \alpha_j(t)(\hat{\lambda}_2(t,\epsilon))_j - \frac{1}{2}\alpha_j^2(t)R_{jj}(t,\epsilon) & \text{if } \alpha_j(t) > (R(t,\epsilon)^{-1}\hat{\lambda}_2(t,\epsilon))_j, \\ \beta_j(t)(\hat{\lambda}_2(t,\epsilon))_j - \frac{1}{2}\beta_j^2(t)R_{jj}(t,\epsilon) & \text{if } \beta_j(t) < (R(t,\epsilon)^{-1}\hat{\lambda}_2(t,\epsilon))_j, \end{cases}$$

for $j = 1, \ldots k$. Finally, solving $f_3^*(\hat{\lambda}_3)$ and $f_4^*(\hat{\lambda}_4)$ we obtain,

$$f_3^*(\hat{\lambda}_3) = \hat{\lambda}_3(\epsilon)^T z_0, \tag{86}$$

$$f_4^*(\hat{\lambda}_4) = \frac{1}{2}\hat{\lambda}_4(\epsilon)^T\pi(\epsilon)^{-1}\hat{\lambda}_4(\epsilon). \tag{87}$$

Substituting (83), (85), (86), and (87) into (81), the dual problem becomes

$$\begin{cases} \underset{(\hat{\rho},\hat{\lambda}_2,\hat{\lambda}_3,\hat{\lambda}_4)}{\text{maximise}} \quad J^{\mathbf{D}}(\hat{\rho},\hat{\lambda}_2,\hat{\lambda}_3,\hat{\lambda}_4,\epsilon), \\ \\ \text{subject to } \; (\hat{\rho},\hat{\lambda}_2,\hat{\lambda}_3,\hat{\lambda}_4) \in \mathbf{\Sigma}_1, \end{cases}$$

where

$$J^{\mathbf{D}} = \int_0^1 -\frac{1}{2}\hat{\rho}(t,\epsilon)^T Q(t,\epsilon)^{-1}\hat{\rho}(t,\epsilon) - \theta(\hat{\lambda}_2)dt - \hat{\lambda}_3(\epsilon)^T z_0 - \frac{1}{2}\hat{\lambda}_4(\epsilon)^T\pi(\epsilon)^{-1}\hat{\lambda}_4(\epsilon).$$

The following lemma will allow us to prove Theorem 1. The derivation below is similar to that in [7] but modified to allow for the singularly perturbed dynamics.

LEMMA 7. *The subspace* $\mathbf{\Sigma}_1$*, orthogonal to* $\mathbf{\Sigma}$*, is given by*

$$\mathbf{\Sigma}_1 = \Big\{(\hat\rho, \hat\lambda_2, \hat\lambda_3, \hat\lambda_4) : \hat\rho(t,\epsilon) = I^{\epsilon}\frac{d\hat\gamma}{dt} + A(t,\epsilon)^T\hat\gamma(t,\epsilon),\ \ \hat\lambda_2(t,\epsilon) = b(t,\epsilon)^T I^{\epsilon}\hat\gamma(t,\epsilon),$$
$$\hat\lambda_3(\epsilon) = I^{\epsilon}\hat\gamma(0,\epsilon),\ \ \hat\lambda_4(\epsilon) = -I^{\epsilon}\hat\gamma(1,\epsilon),\ \ t\in[0,1],\ \ \epsilon\in(0,\epsilon^*]\Big\},$$

*where*

$$\hat\gamma(t,\epsilon) = -I^{\frac{1}{\epsilon}}\bigg(\int_t^1 \Phi_{I^{\frac{1}{\epsilon}}A}(s,t)^T\hat\rho(s,\epsilon)dt - \Phi_{I^{\frac{1}{\epsilon}}A}(1,t)^T I^{\epsilon}\hat\gamma(1)\bigg).$$

*and* $\Phi_{I^{\frac{1}{\epsilon}}A}$ *is the resolvent matrix of the differential equations in* $\mathbf{P}$*.*

*Proof.* For simplicity we omit the dependence on $\epsilon$ for all variables and we omit the $\hat{}$ notation. Let $(\rho,\lambda_2,\lambda_3,\lambda_4)\in\mathbf{\Sigma}_1$. As $\mathbf{\Sigma}_1$ is orthogonal to $\mathbf{\Sigma}$

$$\lambda_4^T s_1 + \lambda_3^T s_0 + \int_0^1 \rho(t)^T z(t) + \lambda_2(t)^T u(t)dt = 0. \tag{88}$$

The solution to the differential equations in $\mathbf{P}$ for time $t$ and time $t=1$ is given by

$$\begin{aligned} z(t) &= \Phi_{I^{\frac{1}{\epsilon}}A}(t,0)s_0 + \int_0^t \Phi_{I^{\frac{1}{\epsilon}}A}(t,s)b(s)u(s)ds,\\ s_1 = z(1) &= \Phi_{I^{\frac{1}{\epsilon}}A}(1,0)s_0 + \int_0^1 \Phi_{I^{\frac{1}{\epsilon}}A}(1,s)b(s)u(s)ds, \end{aligned} \tag{89}$$

respectively. Note that the equations in (89) are integrable because the resolvent matrix satisfies the bounds in (71). Substituting (89) into (88) yields

$$\begin{aligned} &\lambda_4^T\bigg(\Phi_{I^{\frac{1}{\epsilon}}A}(1,0)s_0 + \int_0^1 \Phi_{I^{\frac{1}{\epsilon}}A}(1,s)b(s)u(s)ds\bigg) + \lambda_3^T s_0\\ &+\int_0^1\bigg(\rho(t)^T\bigg(\Phi_{I^{\frac{1}{\epsilon}}A}(t,0)s_0 + \int_0^t \Phi_{I^{\frac{1}{\epsilon}}A}(t,s)b(s)u(s)ds\bigg) + \lambda_2(t)^T u(t)\bigg)dt = 0. \end{aligned}$$

Changing the order of integration

$$\begin{aligned} \lambda_4^T\Phi_{I^{\frac{1}{\epsilon}}A}(1,0)s_0 + \lambda_4^T\int_0^1 \Phi_{I^{\frac{1}{\epsilon}}A}(1,s)b(s)u(s)ds + \lambda_3^T s_0 + \int_0^1 \rho(s)^T\Phi_{I^{\frac{1}{\epsilon}}A}(s,0)s_0 ds\\ +\int_0^1\bigg(\int_s^1 \rho(t)^T\Phi_{I^{\frac{1}{\epsilon}}A}(t,s)dt\bigg)b(s)u(s)ds + \int_0^1 \lambda_2^T(s)u(s)ds = 0. \end{aligned}$$

After some manipulation of terms, we obtain

$$\begin{aligned} &\bigg(\lambda_4^T\Phi_{I^{\frac{1}{\epsilon}}A}(1,0) + \lambda_3^T + \int_0^1 \rho(s)^T\Phi_{I^{\frac{1}{\epsilon}}A}(s,0)ds\bigg)s_0\\ &+\int_0^1\bigg(\lambda_2(s)^T + \bigg(\int_s^1 \rho(t)^T\Phi_{I^{\frac{1}{\epsilon}}A}(t,s)dt + \lambda_4^T\Phi_{I^{\frac{1}{\epsilon}}A}(1,s)\bigg)b(s)\bigg)u(s)ds = 0. \end{aligned}$$

Since $s_0$ and $u$ are arbitrary we obtain the following equations

$$\lambda_4^T\Phi_{I^{\frac{1}{\epsilon}}A}(1,0) + \lambda_3^T + \int_0^1 \rho(s)^T\Phi_{I^{\frac{1}{\epsilon}}A}(s,0)ds = 0, \tag{90}$$

$$\lambda_2(s) + b(s)^T\bigg(\int_s^1 \Phi_{I^{\frac{1}{\epsilon}}A}(t,s)^T\rho(t)dt + \Phi_{I^{\frac{1}{\epsilon}}A}(1,s)^T\lambda_4\bigg) = 0. \tag{91}$$

Define

$$\gamma(s) = -I^{\frac{1}{\epsilon}}\left(\int_s^1 \Phi_{I^{\frac{1}{\epsilon}}A}(t,s)^T\rho(t)dt + \Phi_{I^{\frac{1}{\epsilon}}A}(1,s)^T\lambda_4\right). \tag{92}$$

Substituting (92) into (91)

$$\lambda_2(s) = b^T(s)I^\epsilon\gamma(s). \tag{93}$$

Setting $s = 1$ in (92) we get

$$\gamma(1) = -I^{\frac{1}{\epsilon}}\lambda_4. \tag{94}$$

Note that our expression for $\gamma(s)$ in (92) is the general solution to

$$I^\epsilon\dot{\gamma} = -A(s)^T\gamma(s) + \rho(s).$$

Rearranging gives us

$$\rho(s) = I^\epsilon\dot{\gamma}(s) + A(s)^T\gamma(s). \tag{95}$$

Substituting (95) into (90),

$$\lambda_3^T = -\int_0^1 \dot{\gamma}(t)^T I^\epsilon\Phi_{I^{\frac{1}{\epsilon}}A}(t,0)dt - \int_0^1 \gamma(t)^T A(t)\Phi_{I^{\frac{1}{\epsilon}}A}(t,0)dt - \lambda_4^T\Phi_{I^{\frac{1}{\epsilon}}A}(1,0).$$

Integrating by parts gives us

$$\begin{aligned}\lambda_3^T &= -\gamma(t)^T I^\epsilon\Phi_{I^{\frac{1}{\epsilon}}A}(t,0)|_0^1 + \int_0^1 \gamma(t)^T A(t)\Phi_{I^{\frac{1}{\epsilon}}A}(t,0)dt \\ &\quad - \int_0^1 \gamma(t)^T A(s)\Phi_{I^{\frac{1}{\epsilon}}A}(t,0)dt - \lambda_4^T\Phi_{I^{\frac{1}{\epsilon}}A}(1,0), \\ &= -\gamma(1)^T I^\epsilon\Phi_{I^{\frac{1}{\epsilon}}A}(1,0) + \gamma(0)^T I^\epsilon - \lambda_4^T\Phi_{I^{\frac{1}{\epsilon}}A}(1,0).\end{aligned} \tag{96}$$

Substituting (94) into (96)

$$\lambda_3^T = \gamma(0)^T I^\epsilon. \tag{97}$$

From (93), (94), (95), (97), we obtain the equations in $\mathbf{\Sigma}_1$. □

**6. Numerical Experiment.** We carry out two numerical experiments in this section. The first example is a well known aerospace engineering SPOC problem taken from NASA's Digitial Fly-by-Wire (DFW) aircraft program [12] and adapted in [16]. For this problem, we explicitly state the primal, dual and reduced problems. We choose to run numerical experiments on the model developed in [12] and [16] rather than on a more complex model in order to provide a simple example with intuition. The results of this paper, however, are applicable to larger problems where greater computational savings will be achieved. The second example considers a set of 50 randomly generated SPOC problems of the form $\mathbf{P}$ evaluated at various values of $\epsilon$. We compute the average computational time for the primal and dual SPOC problems and for the upper and lower bounds. In addition, we plot the solution to the primal, dual and reduced problems for one particular randomly generated problem and show that the upper and lower bounds converge to the reduced problem as $\epsilon \to 0$. Our goal in the examples is to demonstrate both the tightness of the error bound and the computational time saved when solving for the bounds instead of for the solution of the SPOC problem.

For the numerical optimisation of the primal, dual and reduced problems, we employ the state of the art optimisation software GPOPS-II [25]. We used the default parameters for the computation with the following exceptions: **tolerance**=$10^{-7}$and **derivatives.supplier** = 'adigator'. The approximate solutions were computed in Matlab where the ODE solver ode15s was employed with the following values: **RelTol**= $10^{-5}$, **AbsTol**=$10^{-7}$, **timesteps** = 100000 for $\epsilon < 0.1$ and **timesteps** = 100 for $\epsilon \geq 0.1$. All integration in in Matlab was computed with the same tolerances.

As the primal and dual problems and their reduced counterparts are solved in GPOPS-II and the upper and lower bounds are computed using an ODE solver in Matlab, the upper and lower bounds must be considered to have some numerical error with respect to the GPOPS-II solution. In particular, depending on the error tolerance specified for ode15s, it is possible for the upper bound to be smaller than the primal solution and the lower bound to be larger. Both the upper and lower bound, however, consistently converge to the reduced solution as $\epsilon \to 0$. When the bounds are calculated using GPOPS-II with the same error tolerance, there is no saving in computational time but the upper and lower bounds are guaranteed to hold and converge to the reduced solution. The CPU times for the primal and dual problems were obtained using the internal GPOPS-II timer and the CPU times for the bounds were obtained using **tic-toc**.

**6.1. Example 1.** The first example concerns the DFW program developed by NASA. The DFW interprets the pilots' flight path input signals and picks the flight control that will achieve this path while accounting for both the aircraft dynamics and the various sensors relating to aircraft performance. We have adapted the SPOC problem in [12] and [16] to fit the form in $\mathbf{P}$. The problem is described as follows

$$
\begin{aligned}
\underset{\hat{u}_1,\hat{u}_2\in U}{\text{minimise}}\quad & \frac{1}{2}\int_0^{60}\Big(\hat{x}_1^2+\hat{x}_2^2+\hat{y}_1^2+\hat{y}_2^2+\hat{u}_1^2+\hat{u}_2^2\Big)dt+\frac{1}{2}\Big(\hat{x}_1(1)^2+\hat{x}_2(1)^2+\epsilon\hat{y}_1(1)^2+\epsilon\hat{y}_2(1)^2\Big)\\
\text{subject to}\quad & \dot{\hat{x}}_1=-0.015\hat{x}_1-0.0805\hat{x}_2-0.00116666\hat{y}_1-0.00009\hat{u}_1+0.02225\hat{u}_2,\\
& \dot{\hat{x}}_2=0.0333333\hat{y}_2,\\
& \epsilon\dot{\hat{y}}_1=-0.076\hat{x}_1-0.028\hat{y}_1+0.0333333y_2-\epsilon0.11\hat{u}_1,\\
& \epsilon\dot{\hat{y}}_2=0.02\hat{x}_1-0.16\hat{y}_1-0.0163333\hat{y}_2-\epsilon8.7\hat{u}_1,\\
& U=[0,1]\times[0,1],\\
& \hat{x}_1(0)=1.55,\ \ \hat{x}_2(0)=0.2,\ \ \hat{y}_1(0)=9,\ \ \hat{y}_2(0)=15,
\end{aligned}
$$

where $\epsilon = 1/30$. Here, $\hat{x}_1$ is the velocity (ft/s), $\hat{x}_2$ is the pitch attitude (rad), $\hat{y}_1$ is the angle of attack (rad), and $\hat{y}_2$ is the pitch rate (rad/s). The control $\hat{u}_1$ is the elevator input which changes the pitch and the control $\hat{u}_2$ is the throttle position. The above variables represent the various quantities relative to a pre-defined equilibrium position of the aircraft. We wish to find the optimal control over a time horizon of 60 sec that returns the aircraft to the equilibrium position. The reduced problem, $\bar{\mathbf{P}}$, becomes

$$\begin{aligned}
\underset{\tilde{u}_1,\tilde{u}_2\in U}{\text{minimise}} \quad & \frac{1}{2}\int_0^{60}\left(5.8351\tilde{x}_1^2+\tilde{x}_2^2+\tilde{u}_1^2+\tilde{u}_2^2\right)dt+\frac{1}{2}\left(\tilde{x}_1(1)^2+\tilde{x}_2(1)^2\right)\\
\text{subject to} \quad & \dot{\tilde{x}}_1=-0.01488\tilde{x}_1-0.0805\tilde{x}_2-0.00009\tilde{u}_1+0.02225\tilde{u}_2,\\
& \dot{\tilde{x}}_2=0.07322\tilde{x}_1,\\
& U=[0,1]\times[0,1],\\
& \tilde{x}_1(0)=1.55, \quad \tilde{x}_2(0)=0.2.
\end{aligned}$$

Let $\hat{\rho}^T=\left[\hat{\rho}_1,\hat{\rho}_2,\hat{\rho}_3,\hat{\rho}_4\right]$, $\hat{\gamma}^T=\left[\hat{\eta}_1,\hat{\eta}_2,\hat{\nu}_1,\hat{\nu}_2\right]$ and $z_0^T=\left[1.55,0.2,9,15\right]$. The dual problem $\mathbf{D}$ becomes

$$\begin{aligned}
\underset{\hat{\rho},\hat{\gamma}}{\text{maximise}} \quad & -\frac{1}{2}\int_0^{60}\left(\hat{\rho}^T\hat{\rho}+\sum_{i=1}^{2}\theta_i(\hat{\gamma})\right)dt-\frac{1}{2}I^\epsilon\hat{\gamma}^T\hat{\gamma}-I^\epsilon\hat{\gamma}(0)^Tz_0\\
\text{subject to} \quad & \dot{\hat{\eta}}_1=0.015\hat{\eta}_1+0.076\hat{\nu}_1-0.02\hat{\nu}_2+\hat{\rho}_1,\\
& \dot{\hat{\eta}}_2=0.0805\hat{\eta}_1+\hat{\rho}_2,\\
& \epsilon\dot{\hat{\nu}}_1=0.00116666\hat{\eta}_1+0.028\hat{\nu}_1+0.16\hat{\nu}_2+\hat{\rho}_3,\\
& \epsilon\dot{\hat{\nu}}_2=-0.0333333\hat{\eta}_2-0.0333333\hat{\nu}_1+0.0163333\hat{\nu}_2+\hat{\rho}_4,
\end{aligned}$$

where

$$\theta_1(\hat{\gamma})=\begin{cases}0 & \text{if } -0.00009\hat{\eta}_1-\epsilon0.11\hat{\nu}_1-\epsilon8.7\hat{\nu}_2<0,\\ \frac{1}{2}(-0.00009\hat{\eta}_1-\epsilon0.11\hat{\nu}_1-\epsilon8.7\hat{\nu}_2)^2 & \text{if } 0\le-0.00009\hat{\eta}_1-\epsilon0.11\hat{\nu}_1-\epsilon8.7\hat{\nu}_2\le1,\\ -0.00009\hat{\eta}_1-\epsilon0.11\hat{\nu}_1-\epsilon8.7\hat{\nu}_2-\frac{1}{2} & \text{if } -0.00009\hat{\eta}_1-\epsilon0.11\hat{\nu}_1-\epsilon8.7\hat{\nu}_2>1,\end{cases}$$

$$\theta_2(\hat{\gamma})=\begin{cases}0 & \text{if } 0.02225\hat{\eta}_1<0,\\ \frac{1}{2}0.02225\hat{\eta}_1^2 & \text{if } 0\le0.02225\hat{\eta}_1\le1,\\ 0.02225\hat{\eta}_1-\frac{1}{2} & \text{if } 0.02225\hat{\eta}_1>1.\end{cases}$$

The numerical performance of the primal, dual, reduced problems along with the upper and lower bounds is summarised in Table 1 for $\epsilon=1/30$. We report the computational time as $\epsilon\to0$ in Table 2.

TABLE 1
*Example 1 - CPU times of* $\mathbf{P}$, $\mathbf{D}$, $\bar{\mathbf{P}}$, *and bounds for* $\epsilon=1/30$.

| CPU Time [Sec] | | | |
|---|---|---|---|
| **P** | **D** | $\bar{\mathbf{P}}$ | Upper Bd and Lower Bd |
| 2.95 | 4.56 | 2.83 | 4.13 |

The solutions to the primal, dual and reduced problems and the upper and lower bounds are plotted as $\epsilon\to0$ in Fig. 2. The difference between the upper and lower bounds is plotted in Fig. 3. It is clear from the graph that one would consider the reduced solution a good approximation to the solution to the original problem when $\epsilon>6*10^{-3}$. Furthermore, we see that the primal and dual solutions give the same objective function values and that the upper and lower bounds are converging to the reduced solution as $\epsilon\to0$ within some numerical error. Calculating the constant $C$ in (15) yields

$$C=\frac{1}{140.5011}\left(\frac{140.6390-140.4621}{0.00001}\right)=125.9065$$

TABLE 2
*Example 1 - CPU times of* **P** *and* **D** *and bounds as* $\epsilon \to 0$.

| $\epsilon$ | CPU Time [Sec] | | | Gain **P** vs. bounds | Gain **D** vs. bounds |
|---|---|---|---|---|---|
| | **P** | **D** | Upper Bd and Lower Bd | | |
| 0.01 | 7.23 | 8.21 | 4.96 | 1.46 | 1.65 |
| 0.001 | 45.22 | 16.87 | 6.28 | 7.2 | 2.69 |
| 0.0001 | 15.95 | 50.11 | 4.7 | 3.39 | 10.66 |
| 0.00001 | 41.75 | 132.94 | 4.92 | 8.49 | 27.05 |

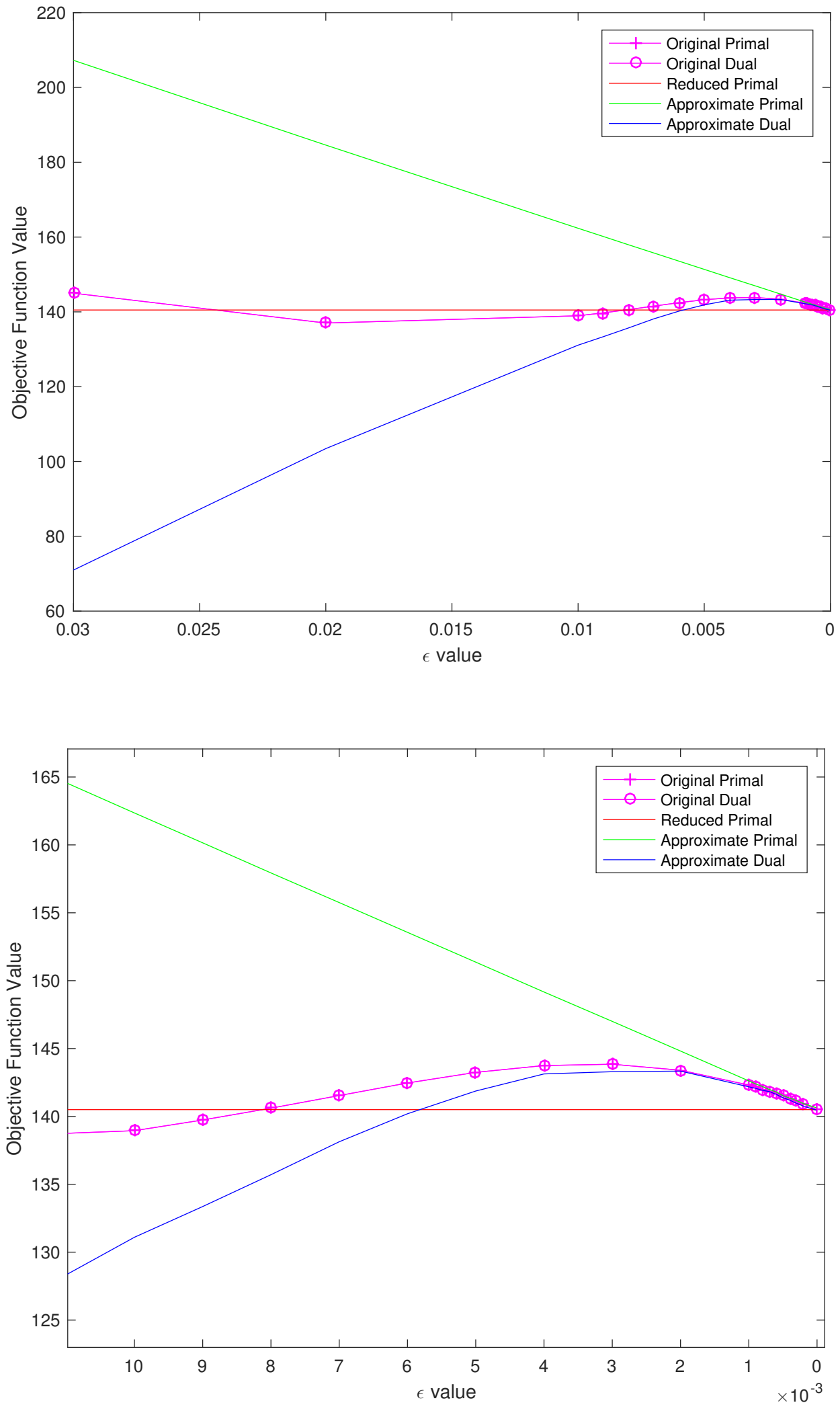

FIG. 2. *Example 1 - Objective Function Values for Original and Reduced problems and bounds.*

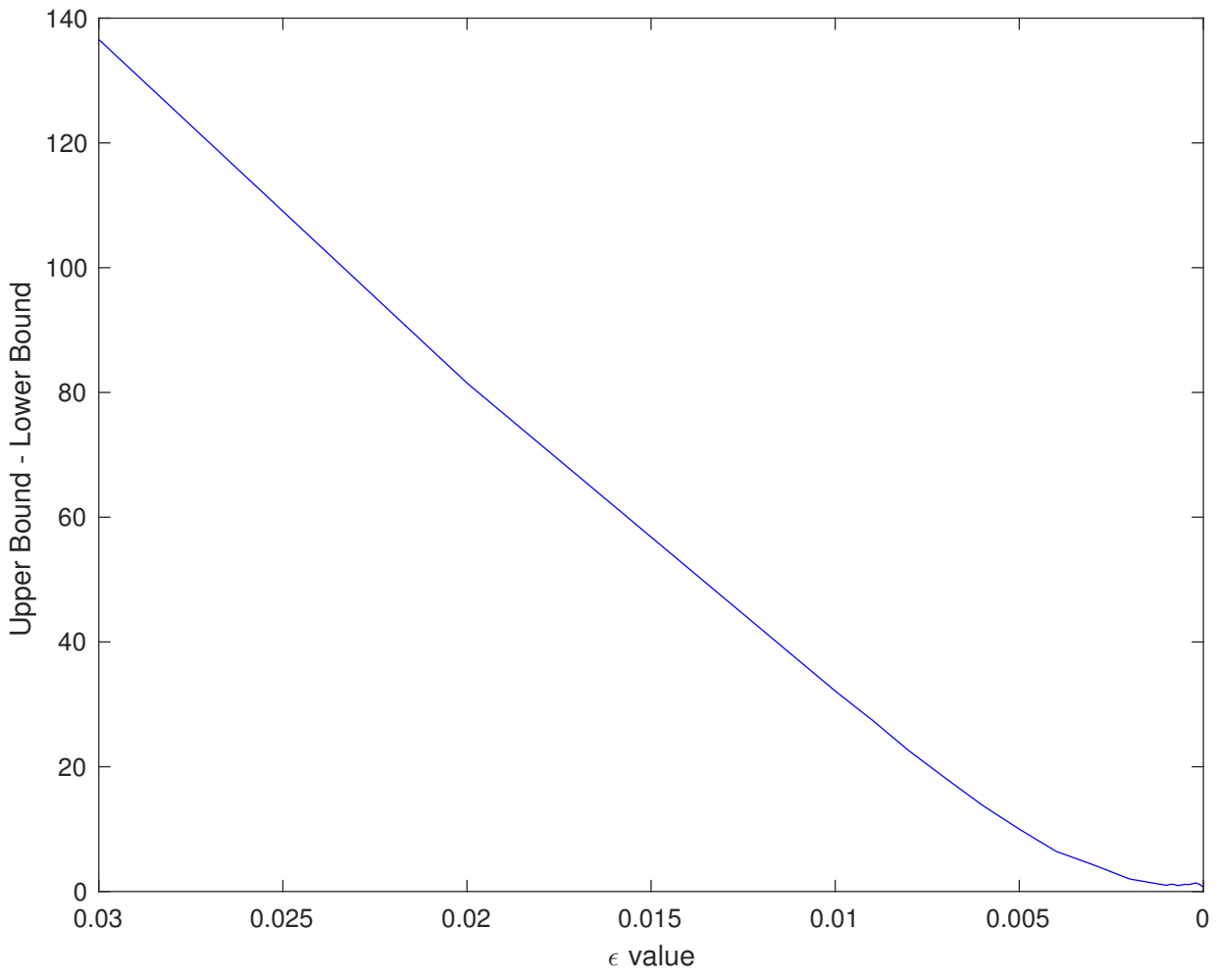

FIG. 3. *Example 1 - Difference between the upper and lower bounds as $\epsilon \to 0$.*

**6.2. Example 2.** We consider problems of the form **P** with the values $A$, $b$, $Q$, $R$, $\pi$, $\alpha$, $\beta$, $z_0$ randomly generated and such that the corresponding reduced problems are well conditioned. For all problems, we take $\hat{z_1} \in W^{1,2}([0,0.5];\mathbb{R}^4)$, $\hat{z_2} \in W^{1,2}([0,0.5];\mathbb{R}^6)$ and $\hat{u} \in L^2([0,0.5];\mathbb{R}^3)$. We run 50 examples and report the average computational time for the primal problem, the dual problem and the upper and lower bounds for various values of $\epsilon$ in Table 3.

TABLE 3
*Example 2 - Average CPU times of* **P** *and* **D** *and bounds*

| $\epsilon$ value | CPU Time [Sec] | | | Gain |
|---|---|---|---|---|
| | **P** | **D** | Upper Bd and Lower Bd | **D** vs bounds |
| 1 | 0.88 | 0.84 | 2.16 | 0.39 |
| 0.1 | 2.77 | 1.87 | 2.02 | 0.93 |
| 0.01 | 24.11 | 15.37 | 4.18 | 3.68 |
| 0.001 | 106.15 | 25.71 | 6.64 | 3.87 |
| 0.0001 | 720.01 | 40.98 | 6.74 | 6.08 |
| 0.00001 | 1160.13 | 141.45 | 7.13 | 19.84 |

As $\epsilon \to 0$, there is a significant increase in computational time when solving the primal rather than the dual and in some instances, the solver was not able to produce a solution to the primal problem. Similar results for the unperturbed primal and dual case were produced in [7]. The computational time for the upper and lower bounds for small $\epsilon$, however, is consistently less than the computational time taken to solve either the dual or the primal problem.

We plot the objective function values for one particular randomly generated problem. The values used to generate this problem are recorded in the Appendix. Fig. 4 shows the objective function values for the primal, dual and reduced problems along with the bounds as $\epsilon \to 0$. The difference between the upper and lower bound is plotted in Fig. 5. It is clear that the upper and lower bounds hold for arbitrary $\epsilon$ and converge to the reduced solution as $\epsilon \to 0$. Futhermore, for this problem, we would consider the reduced solution to provide a poor approximation to the primal solution on the entire interval. Calculating the constant $C$ in (15) yields

$$C = \frac{1}{346.2091}\left(\frac{346.2090 - 344.0990}{0.00001}\right) = 10.3868.$$

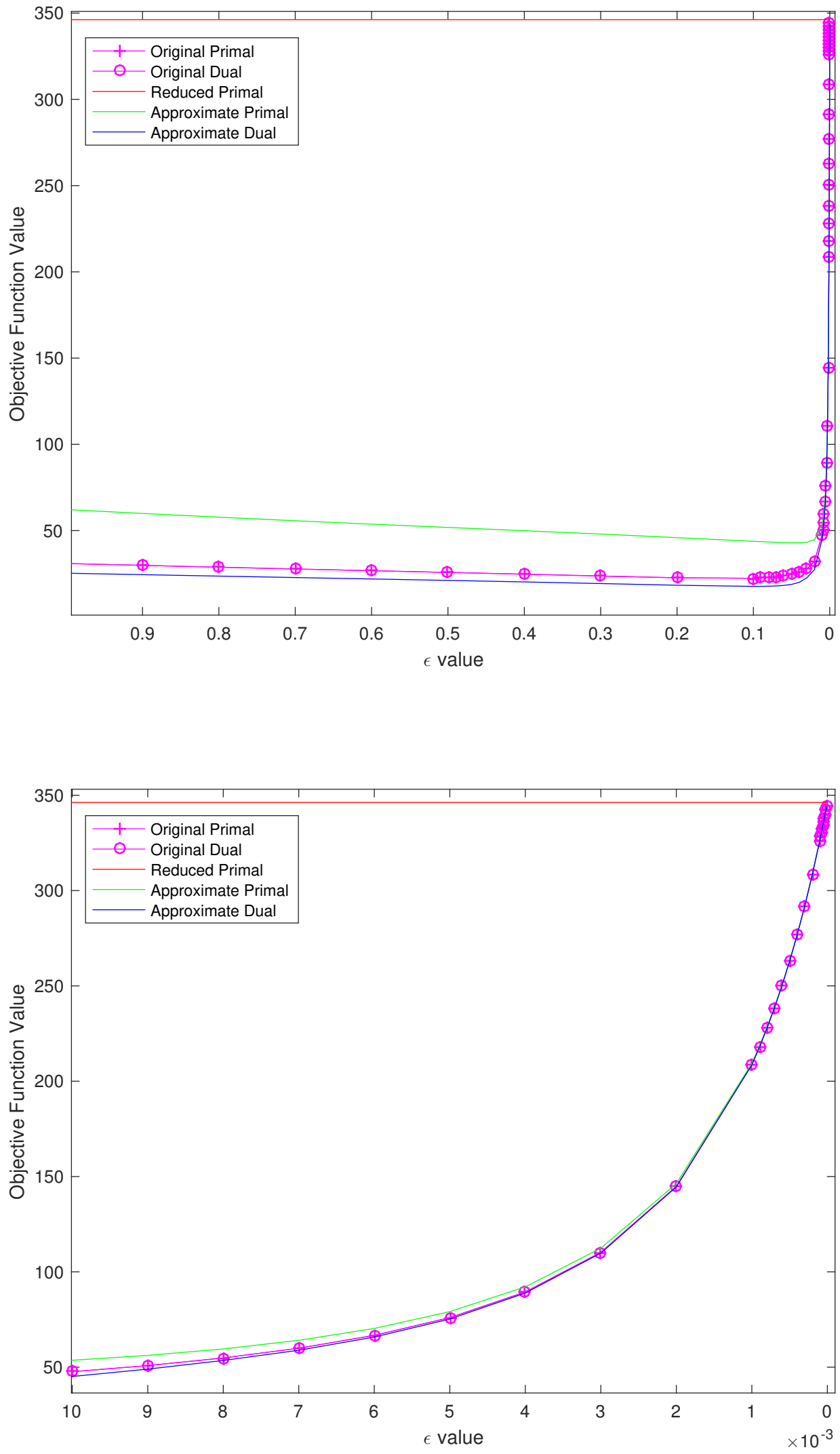

FIG. 4. *Example 2 - Objective function values for the Original and Reduced problems and bounds.*

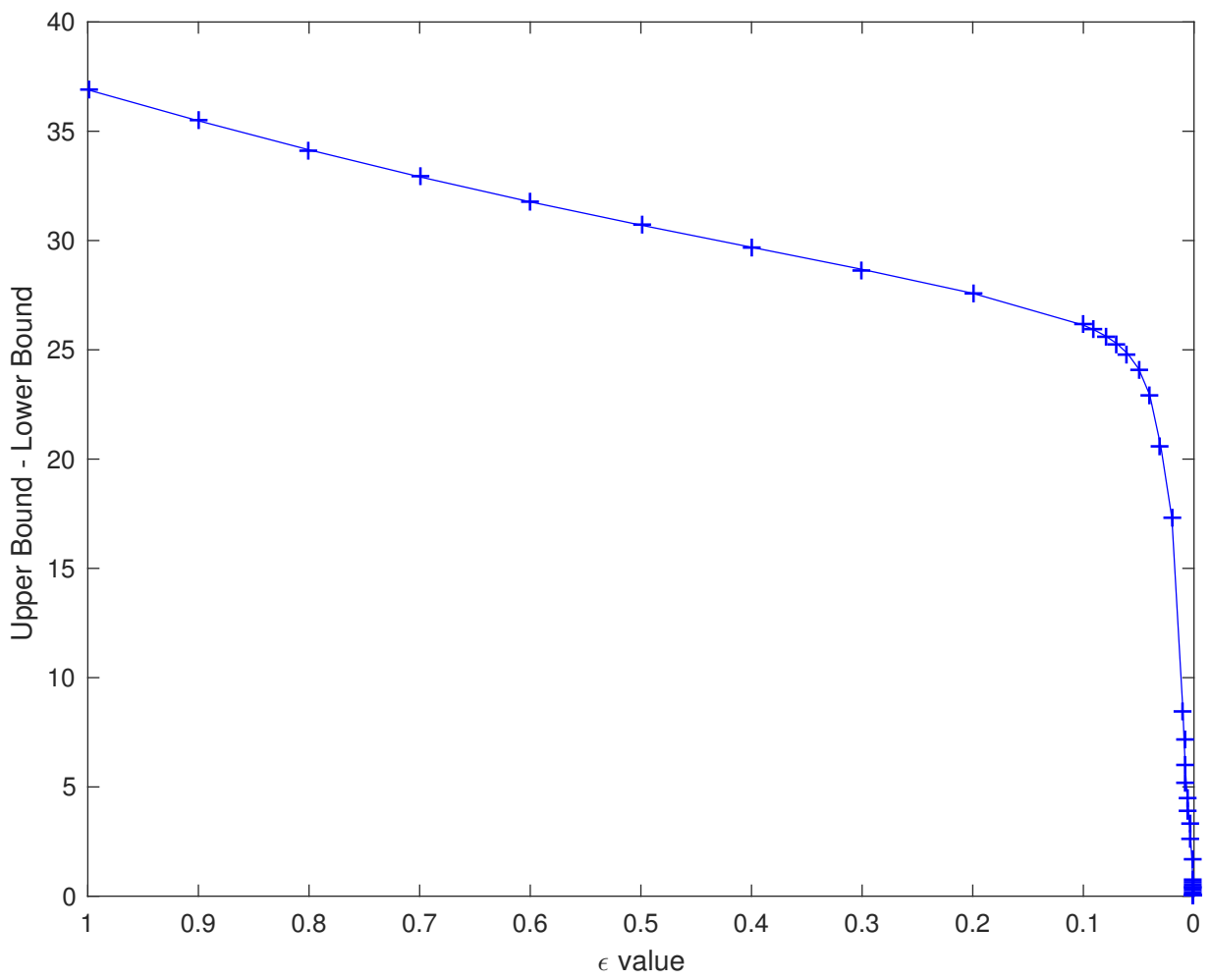

FIG. 5. *Example 2 - Difference between the upper and lower bound as* $\epsilon \to 0$.

**7. Conclusion.** We have developed a methodology to compute upper and lower bounds on the solution of a SPOC problem such that the ensuing error term is bounded by $C\epsilon$ for some known constant $C$. This methodology rests on our construction of a dual SPOC problem whose optimal solution converges to that of a reduced SPOC problem and on our extension of the strong duality theorem from the result in [7] to SPOC problems. Moreover, we have developed a quantifiable criterion for establishing how good of an approximation the reduced solution provides. In addition, our numerical experiments show that the computational time saved when evaluating the upper and lower bounds, instead of the original SPOC problem, can be significant even for small problems and they demonstrate the convergence of both bounds to the optimal solution of the SPOC problem as $\epsilon \to 0$.

**Appendix.** The following values were randomly generated and were used in Example 2 to show the objective function values of the original and reduced problems and the associated upper and lower bound.

$$\alpha = \begin{bmatrix} 2.775 & 3.441 & 3.485 \end{bmatrix},$$
$$\beta = \begin{bmatrix} 5.105 & 6.470 & 5.981 \end{bmatrix},$$
$$z_0 = \begin{bmatrix} -0.437 & 0.940 & -0.8180 & 0.138 & 0.7070 & 0.054 & 0.965 & 0.665 & 0.130 & 0.292 \end{bmatrix},$$
$$A_{11} = \begin{bmatrix} 0.905 & 0.289 & -0.209 & -0.226 \\ -0.696 & 0.973 & -0.194 & -0.706 \\ -0.806 & -0.622 & -0.868 & 0.198 \\ -0.988 & -0.407 & 0.469 & 0.026 \end{bmatrix},$$

$$A_{12} = \begin{bmatrix} -0.210 & -0.774 & -0.687 & 0.208 & -0.131 & -0.693 \\ -0.376 & -0.095 & 0.043 & 0.817 & -0.500 & -0.406 \\ -0.256 & 0.860 & -0.944 & 0.999 & 0.838 & -0.867 \\ -0.305 & -0.627 & 0.136 & -0.103 & 0.920 & -0.764 \end{bmatrix},$$

$$A_{21} = \begin{bmatrix} 0.566 & -0.121 & -0.784 & -0.649 \\ -0.344 & -0.954 & 0.855 & -0.393 \\ 0.1756 & -0.742 & -0.837 & 0.622 \\ 0.792 & -0.581 & -0.174 & -0.124 \\ 0.558 & -0.481 & -0.363 & -0.147 \\ 0.226 & 0.146 & -0.447 & -0.533 \end{bmatrix},$$

$$A_{22} = \begin{bmatrix} -2.380 & -0.773 & 0.674 & -1.261 & -0.144 & 0.280 \\ -0.773 & -1.958 & 1.163 & -0.003 & -0.603 & 1.676 \\ 0.674 & 1.163 & -2.388 & 1.099 & -0.506 & -0.614 \\ -1.261 & -0.003 & 1.099 & -3.085 & 0.777 & -0.711 \\ -0.144 & -0.603 & -0.506 & 0.777 & -1.045 & 1.030 \\ 0.280 & 1.676 & -0.614 & -0.711 & 1.030 & -1.943 \end{bmatrix},$$

$$b_1 = \begin{bmatrix} -0.467 & -0.302 & -0.886 \\ 0.768 & -0.261 & -0.374 \\ 0.910 & -0.077 & 0.873 \\ -0.936 & -0.472 & 0.140 \end{bmatrix},$$

$$b_2 = \begin{bmatrix} 0.945 & 0.166 & -0.266 \\ 0.0198 & -0.917 & 0.047 \\ 0.890 & -0.037 & 0.390 \\ 0.288 & 0.111 & -0.366 \\ 0.066 & 0.164 & 0.279 \\ -0.134 & 0.476 & 0.302 \end{bmatrix}$$

$$R = \begin{bmatrix} 0.124 & 0 & 0 \\ 0 & 0.422 & 0 \\ 0 & 0 & 0.107 \end{bmatrix},$$

$$\pi_{11} = \begin{bmatrix} 1.661 & -0.492 & 0.467 & -0.571 \\ -0.492 & 1.403 & 0.509 & 0.689 \\ 0.467 & 0.509 & 0.908 & -0.049 \\ -0.571 & 0.689 & -0.049 & 1.42 \end{bmatrix},$$

$$\pi_{22} = \begin{bmatrix} 2.695 & 1.886 & 0.342 & 1.018 & -1.920 & -0.792 \\ 1.886 & 2.136 & -0.724 & 1.051 & -1.140 & -0.525 \\ 0.342 & -0.723 & 2.109 & -0.488 & -0.597 & 0.973 \\ 1.018 & 1.051 & -0.488 & 1.064 & -0.492 & -1.093 \\ -1.920 & -1.140 & -0.597 & -0.492 & 2.324 & 0.484 \\ -0.792 & -0.525 & 0.973 & -1.093 & 0.484 & 2.289 \end{bmatrix},$$

$$Q = \begin{bmatrix} 4.334 & 2.848 & 0.617 & -0.611 & -0.104 & 0.730 & 0.679 & -1.793 & 1.870 & 0.445 \\ 2.848 & 3.008 & 0.714 & -1.246 & 0.822 & 0.157 & 0.279 & -0.868 & 1.316 & -0.550 \\ 0.617 & 0.714 & 5.423 & 0.292 & 2.128 & 1.009 & 0.371 & -2.352 & 2.926 & -1.511 \\ -0.611 & -1.246 & 0.292 & 3.898 & -0.026 & 0.589 & -1.184 & -2.039 & -0.119 & 0.574 \\ -0.104 & 0.822 & 2.128 & -0.026 & 2.643 & 0.029 & 0.227 & -1.541 & 1.125 & -0.676 \\ 0.730 & 0.157 & 1.009 & 0.589 & 0.029 & 4.287 & 0.652 & -0.367 & 1.341 & 1.589 \\ 0.679 & 0.279 & 0.371 & -1.184 & 0.227 & 0.652 & 2.621 & -0.213 & 1.204 & 0.683 \\ -1.793 & -0.868 & -2.352 & -2.039 & -1.541 & -0.367 & -0.213 & 3.907 & -1.814 & -0.350 \\ 1.870 & 1.316 & 2.926 & -0.120 & 1.125 & 1.341 & 1.204 & -1.814 & 3.040 & -0.620 \\ 0.445 & -0.550 & -1.511 & 0.574 & -0.676 & 1.589 & 0.683 & -0.350 & -0.620 & 3.206 \end{bmatrix}.$$

**References.**